# Nonexchangeable logit dynamics with cost constraints: different viewpoints and numerical analysis

Hidekazu Yoshioka[1,*], Tomohiro Tanaka[2]

[1] Graduate School of Advanced Science and Technology, Japan Advanced Institute of Science and Technology, 1-1 Asahidai, Nomi, Ishikawa, Japan

[2] Disaster Prevention Research Institute, Kyoto University, Gokasho, Uji, Kyoto, Japan

[*] Corresponding author: yoshih@jaist.ac.jp, ORCID: 0000-0002-5293-3246

***Abstract***

The logit dynamic, a nonlinear dynamical system, subject to a const constraint with heterogeneous agents is formulated, and its well-posedness is studied from both agent-based (stochastic differential equation: SDE) and probabilistic (Fokker–Planck equation: FPE) viewpoints. The state space of agent actions is a compact domain in a finite-dimensional Euclidean space. The SDE is of the McKean–Vlasov type and is driven by jumps with finite variations, whereas the FPE is a nonlinear integro-differential equation. A key to our mathematical analysis of the FPE is adding a regularization factor into the cost constraint to mitigate the blow-up of the logit function. This property carries over to the McKean–Vlasov SDE. We also present a numerical method based on a naïve finite difference discretization for computing the FPE along with demonstrative computational examples, showing that the regularization method does not critically affect numerical solutions while preventing the breakdown of numerical computation. Finally, we conduct another demonstrative application study in which environmental, energy, and fishery resources intersect.

***Statements & Declarations***

**Competing interests:** The authors have no relevant financial or nonfinancial interests.

**Data availability:** Data will be made available upon reasonable request to the author.

**Funding:** Japan Society for the Promotion of Science (KAKENHI, No. 25K00240, 26K00357)

**Acknowledgments:** This study was supported by Japan Society for the Promotion of Science (KAKENHI, No. 25K00240, 26K00357).

**Declaration of Generative AI and AI-assisted technologies in the writing process:** The author used generative AI to proofread the manuscript.

**Contributions: (Hidekazu Yoshioka)** Conceptualization, Methodology, Software, Writing- Original draft preparation, Visualization, Investigation, Supervision, Writing- Reviewing and Editing, Funding acquisition. **(Tomohiro Tanaka)** Software, Data curation, Visualization, Writing- Reviewing and Editing

## 1. Introduction

### 1.1 Logit dynamics

Evolutionary games are a mathematical framework for analyzing diverse social dynamics [1,2], and logit dynamics are among the major mathematical models used to describe how agents sequentially develop their actions while interacting with each other through utilities. When considering continuous-time dynamics where the agent's action space is a compact set $D$ in a finite-dimensional Euclidean space (e.g., a closed interval), the Fokker–Planck equation (FPE) governing the logit dynamic is given as follows [3] (a more rigorous formulation of each equation in **Section 1** is presented in **Section 2**):

$$\underbrace{\frac{\mathrm{d}\mu_t}{\mathrm{d}t}}_{\text{Temporal change}} = \underbrace{\frac{\exp(\beta U(x,\mu_t))}{\int_D \exp(\beta U(z,\mu_t))\mathrm{d}z}}_{\text{Logit function: inflow}} - \underbrace{\mu_t}_{\text{Linear decay: Outflow}} \left(= \mathbb{L}_\beta(x,\mu_t) - \mu_t\right),\ t>0 \tag{1}$$

subject to some initial condition, where $\mu$ is the probability measure of agent actions and $U$ is the utility. The logit function $\mathbb{L}_\beta$ is a softmax function of utility $U$ whose steepness is controlled by the inverse temperature $\beta > 0$, and $\mathbb{L}_\beta$ is more concentrated at its maxima as $\beta$ increases. Note that discrete-state cases, such as transportation [4,5], cybersecurity [6,7], technological change dynamics [8], cooperative evacuation behavior [9], and sociohydrology [10,11], can be considered if we replace the domain $D$ with a set of finite points with slight modifications. Note also that the logit function $\mathbb{L}_\beta$ can be seen as a probability density function (PDF) on $D$. Logit dynamics have also been studied in the context of machine learning, where they are called best response dynamics or entropic fictitious play [12,13], and their applicability to neural networks has been investigated [14].

The well-posedness (unique existence of solutions) of the logit dynamic (1) and its generalization where the logit function is defined through generalized exponential functions has been extensively studied on the basis of generalized Picard–Lindelöf's theorem [15,16]; see also p.79 in Zeidler [17], which assumes some Lipschitz continuity of utilities with respect to actions and their distributions. Transient solution behaviors of logit dynamics have been reported in the literature, focusing on discrete states [18,19] and problems with delay [20,21].

In the analysis of logit dynamics, studies have also been conducted on the equilibrium solution: at a stationary state ($t = \infty$), the logit dynamic (1) reduces to the fixed-point equation of $\mu_\infty$:

$$\mu_\infty(x) = \mathbb{L}_\beta(x,\mu_\infty), \tag{2}$$

and $\mu_\infty$ approximates a Nash equilibrium when $\beta$ is sufficiently large under certain regularity conditions about $U$, such as certain continuity and concavity (e.g., Theorem 3.2 in Lahkar and Riedel [3]; Corollary 2.4 in Lascu and Majka [12]). In this view, a model with a larger $\beta$ is more accurate in terms of the approximation of the Nash equilibrium. Relationships analogous to (2) also hold true in different mathematical models, such as the Fisher–Rao gradient flows [22,23], which can be seen as replicator dynamics with a logarithmic regularization term. The relationship (2) can also be related to the stationary states of consensus-based optimization models [24,25]. Thus, there are linkages between logit and seemingly different models.

### 1.2 Unresolved issues

We will now describe three unresolved issues concerning logit dynamics.

#### 1.2.1 Exploration cost

On the other hand, the following representation of the logit function provides another interpretation of the inverse temperature $\beta$ as a weighting coefficient in an optimization problem (Section 4 in Lahkar and Riedel [3]):

$$\mathbb{L}_{\beta}(\cdot,\mu_t)=\arg\max_{u}\left\{\underbrace{\int_D U(z,\mu_t)u(z)\mathrm{d}z}_{\text{Average utility}}-\frac{1}{\beta}\underbrace{\int_D\{u(z)\ln u(z)-(u(z)-1)\}\mathrm{d}z}_{\text{Relative entropy: exploration cost}}\right\}, \tag{3}$$

where $u$ in argmax is chosen from the probability densities in $D$. This is an unconstrained maximization problem of the average utility minus exploration cost, where the relative entropy is always nonnegative and becomes larger as $u$ deviates more from a uniform distribution, and the deviation is more strongly penalized for a smaller $\beta$.

In the optimization problem described above, the higher the reverse temperature $\beta$ is, the more diverse strategies the agent can explore. The classical duality argument for unconstrained optimization problems shows that there exists a constant $\varepsilon>0$ depending on $\beta$ such that

$$\mathbb{L}_{\beta}(\cdot,\mu_t)=\arg\max_{u}\left\{\int_D U(z,\mu_t)u(z)\mathrm{d}z:\int_D\{u(z)\ln u(z)-(u(z)-1)\}\mathrm{d}z\leq\varepsilon\right\}. \tag{4}$$

This is a way of expressing the unconstrained optimization problem (3), where the cost, which is called exploration cost in the sequel, is now explicitly constrained. If we understand that $\varepsilon$ is specified first and then the corresponding $\beta$ is found provided that the cost constrained is satisfied with equality, then this $\beta$ also depends on $\mu_t$ through (4). Therefore, the resulting logit dynamic where $\beta$ is formally expressed as $\beta(\mu_t;\varepsilon)$ is more strongly nonlinear than the classical model:

$$\frac{\mathrm{d}\mu_t}{\mathrm{d}t}=\mathbb{L}_{\beta(\mu_t;\varepsilon)}(x,\mu_t)-\mu_t,\ \ t>0. \tag{5}$$

This type of logit dynamics was first studied by Yoshioka [26] under an unusual strong assumption about profiles of utilities $U$ that they are sufficiently far from constant function and are bounded from both below and above in terms of entropy. If $\varepsilon>0$ and $U$ is a constant function, then no such $\beta(\mu_t;\varepsilon)$ exists, and the representation (5) becomes meaningless. This implies that some regularization methods are necessary to formulate meaningful logit dynamics under mild conditions for utilities.

#### 1.2.2 Agent heterogeneity

The logit dynamic (1) assumes that agents are homogeneous and exchangeable. However, when various populations in the real world are considered, the assumption that there are absolutely no differences between individuals is considered extremely unfounded. The evolution of interacting heterogeneous agents has recently been studied, focusing on swarm and oscillator dynamics [27,28,29] and linear-quadratic games [30], with governing equations being different from logit ones, and on "static" (and hence strange) evolutionary game dynamics, where problems involve equilibrium states only [31,32]. While extending classical logit dynamics to models of nonexchangeable agents is theoretically possible, there are no research examples that have expanded the scope to include cases with exploration cost constraints.

#### 1.2.3 McKean–Vlasov dynamics

We consider a complete filtered probability space $\left(\Omega, \mathcal{F}, (\mathcal{F}_t)_{t\geq 0}, \mathbb{P}\right)$ as usual in stochastic calculus, and the expectation is denoted by $\mathbb{E}$. As evidenced by the fact that a logit dynamic can be understood as a (nonlinear) FPE, a stochastic differential equation (SDE) is associated with it (e.g., Chapters 3-4 in Risken [33]). The SDE associated with (1) is formally given as follows, where $X = (X_t)_{t\geq 0}$ is the continuous-time process corresponding to the time-dependent action of a representative agent:

$$\mathrm{d}X_t = \int_D \left(z - X_{t-}\right) M\left(\mathrm{d}t, \mathrm{d}z\right), \tag{6}$$

where $M$ is a Poisson random measure such that $\mathbb{E}\left[M\left(\mathrm{d}t, \mathrm{d}z\right) | \mathcal{F}_{t-}\right] = \mathbb{L}_\beta\left(z, \mu_{t-}\right)\mathrm{d}t\mathrm{d}z$ ($\mathcal{F}_{t-}$ is the left limit of a natural filtration generated by $M$ at time $t$), i.e., the compensator of $M$ is proportional to $\mathbb{L}_\beta$. According to the SDE (6), $X_{t-}$ is updated to $X_t$ drawn from the PDF $\mathbb{L}_\beta\left(z, \mu_{t-}\right)$; formally, we can write $X_t \sim \mathbb{L}_\beta\left(z, \mu_{t-}\right)$ at each jump time. An important observation is that the jump intensity equals 1 and is therefore strictly bounded because $\mathbb{L}_\beta$ is a PDF (e.g., Chapter 1 in Øksendal and Sulem [34]).

When $\mu$ is assumed to be given as a solution to the logit dynamic (1), then the SDE (6) is simply a model driven by a compound Poisson process; in contrast, if $\mu$ is the law of $X$, then the SDE (6) is of the McKean–Vlasov type driven by jump and becomes a nontrivial mathematical object [35]. Furthermore, the problem becomes more complex when the nonexchangeability of agents is assumed. Recently, McKean–Vlasov SDEs for nonexchangeable agent dynamics driven by Brownian motions were studied [36,37], while those driven by jumps have not been explored well.

### 1.3 Purpose and contribution

Considering the three research issues mentioned above (exploration cost, agent heterogeneity, and McKean–Vlasov SDEs), the aims of this study are mathematical and numerical analyses of logit dynamics considering exploration cost constraints based on nonexchangeable agents, as well as a mathematical analysis of the associated McKean–Vlasov SDE. Our contributions to achieving these aims are explained below.

### 1.3.1 Nonexchangeable logit dynamic

Throughout this paper, as discussed later, for utility $U$, we assume Lipschitz continuity for the action $x$ and Lipschitz continuity in a total variation norm for the distribution $\mu$. We also assume that utilities represent the heterogeneity of agents. We show that by introducing a regularization term into the exploration cost constraint, the (solution-dependent) inverse temperature $\beta$ becomes strictly and uniformly bounded from above, with which we can prove the unique existence of solutions to our logit dynamic. A key here is how to design a simple regularization term so that the problem becomes well-posed. The mathematical techniques employed here are based on those used in Yoshioka [26], but the difference here is that this paper does not assume severe assumptions about profiles of $U$. Moreover, we extend the previous result to a nonexchangeable case. This paper imposes abstract assumptions about utilities but not their profiles and therefore covers a wide range of utilities, including potential utilities [38,39] and graphon-based utilities [40,41,42].

### 1.3.2 Nonexchangeable McKean–Vlasov SDE

We analyze the McKean–Vlasov SDE associated with our logit dynamic under a stronger assumption that utility is Lipschitz continuity in a first-order Wasserstein distance, which is a slightly stronger assumption than that with a total variation distance; however, we show that this does not critically affect our problem setting. A major contribution to dealing with the McKean–Vlasov SDE lies in representing the heterogeneity of agents by formally increasing the dimensionality of the Poisson random measure. The Lipschitz continuity in the Wasserstein distance plays a role in obtaining Gronwall-type inequalities. The abstract nature of mathematical analysis for the logit dynamic carries over to the analysis here except for the abovementioned Lipschitz continuity, supporting the consistency of our framework of analysis. Recently, McKean–Vlasov SDEs related to evolutionary games have been mathematically analyzed in the literature [43,44]; however, those related to logit dynamics have not been studied to the best of the authors' knowledge.

### 1.3.3 Computational application

This paper also addresses numerical computation of logit dynamics. Our numerical method uses a naïve finite difference discretization for the logit dynamic and a Newton method for the exploration cost constraint. The regularization method also plays a role here; the existence of the regularization term ensures that the denominator in the Newton iteration, which should have a fixed sign to avoid algorithmic failure, is strictly positive irrespective of the numerical solution values. As demonstrative computational applications, for a test case, we first analyze the influence of regularization parameters on the behavior of numerical solutions from both transient and steady-state perspectives.

Next, by extending the test case above, we tackle a more complex numerical example where the environment, energy, and fishery resources intersect. More specifically, we computationally investigate a fishing problem in a river reach whose streamflow is controlled by an upstream dam. In this application, the action space represents the intensity of fishing and heterogeneity the locations in the river, and the utility is specified according to a suitability index (SI) based on a coupled hydrological–hydraulic model [45]. This kind of combination study between evolutionary game and hydrological–hydraulic models is quite novel.

### 1.4 Structure of this paper

The remainder of this paper is structured as follows. In **Section 2**, we present the logit dynamic with cost constraints under nonexchangeability and the corresponding SDE, along with the main assumptions. In **Section 3**, we describe the results regarding the unique existence of solutions to these equations. The proofs of the propositions in **Section 3** are given in **Appendix A**. In **Section 4**, we present our demonstrable numerical applications by using the logit dynamic. Details of the numerical computation methods are given in **Appendix B**. In **Section 5**, we summarize this study and present its future perspectives.

## 2. Evolutionary game dynamics

### 2.1 Preliminaries

The notations used in the remainder of this paper are explained below. Time $t$ is a nonnegative continuous parameter. The space of agent actions is given by a $n$-dimensional compact domain $D$ with some $n \in \mathbb{N}$. Without any loss of generality, the diameter of $D$ (the maximum Euclidean distance between two points in $D$) is assumed to be 1. The space of agent types is given by the unit interval $\Gamma = [0,1]$ without significant loss of generality (e.g., Yoshioka et al. [31]).

For the logit dynamic, i.e., FPE, we borrow some settings about functional spaces from the literature [31,46]. We set the compact product space $K = D \times \Gamma$. The space of finite signed measures $\mathcal{M}(D)$ is equipped with the total variation norm $\|\mu\|_{\mathrm{TV}} = \sup_{|f| \le 1} \left| \int_D f(x)\mu(\mathrm{d}x) \right|$; here, the supremum is taken with respect to all measurable functions $f : D \to \mathbb{R}$ whose absolute values are bounded by 1. A continuum of signed measures on $D$ parameterized by $y \in \Gamma$ is denoted by $\{\mu(y,\mathrm{d}x)\}_{y \in I}$. The collection of such parameterized measures is denoted by $\mathcal{M}^{(\Gamma)}$, which is equipped with the norm

$$\|\mu\|_{\mathrm{TV}}^{(\Gamma)} = \sup_{y \in I} \|\mu(y,\cdot)\|_{\mathrm{TV}}. \tag{7}$$

For technical reasons, the subspace $\mathcal{M}_2^{(\Gamma)} \subset \mathcal{M}^{(\Gamma)}$ is defined as the collection of all $\mu \in \mathcal{M}^{(\Gamma)}$ such that $\|\mu\|_{\mathrm{TV}}^{(\Gamma)} \le 2$. The space of probability measures on $D$ is denoted by $\mathcal{P}(D)$, and the space of all the continuums of probability measures in $\mathcal{P}(D)$ parameterized by $y \in \Gamma$ is denoted by $\mathcal{P}^{(\Gamma)}$. We also use the first-order Wasserstein distance (e.g., Theorem 1.21 in Chewi et al. [47]): for any $\mu, \nu \in \mathcal{P}(D)$,

$$W_1(\mu,\nu) = \sup_{\phi \in \mathrm{Lip}_1(D)} \int_D \phi(x)\left(\mu(\mathrm{d}x) - \nu(\mathrm{d}x)\right), \tag{8}$$

where $\mathrm{Lip}_1(D)$ is the space of all functions on $D$ such that $|\phi(x_1)-\phi(x_2)| \leq |x_1-x_2|$. By definition, $W_1(\mu,\nu) \leq \|\mu-\nu\|_{\mathrm{TV}}$ for any $\mu,\nu \in \mathcal{P}(D)$. The following upper bound of $W_1$ is convenient (e.g., Proof of Proposition 1.3 in Chewi et al. [47]):

$$W_1(\mu,\nu) \leq \mathbb{E}\left[|X-Z|\right], \tag{9}$$

where $X, Z$ are random variables generated from $\mu, \nu$, respectively.

We also present notations used for studying the McKean–Vlasov SDEs. We consider a complete filtered probability space $\left(\Omega, \mathcal{F}, (\mathcal{F})_{t\geq 0}, \mathbb{P}\right)$ [34]. For each $T>0$, we set the Banach space of the parameterized continuous-time processes valued in $D$ at each time as follows:

$$\mathcal{H}_T = \left\{ \mathbf{X} = \left\{X_t^y\right\}_{t\in[0,T],\, y\in\Gamma} \middle| \|\mathbf{X}\|_{\mathcal{H}_T} < +\infty \right\}, \tag{10}$$

which is equipped with the norm $\|\mathbf{X}\|_{\mathcal{H}_T} = \int_0^1 \sup_{t\in[0,T]} \mathbb{E}\left[\left|X_t^y\right|\right] \mathrm{d}y$.

Finally, for any $a,b \in \mathbb{R}$, we set $\max\{a,b\}$ as $a \vee b$, $\min\{a,b\}$ as $a \wedge b$, and $\max\{a,0\}$ as $(a)_+$. $\mathbb{I}_{\{S\}}$ is the indicator function such that $\mathbb{I}_{\{S\}} = 1$ if $S$ is true and $\mathbb{I}_{\{S\}} = 0$ otherwise.

### 2.2 Utility

A utility is a function $U : D \times \Gamma \times \mathcal{M}^{(\Gamma)} \to \mathbb{R}$. We present assumptions about $U$, **Assumption 1**, which is assumed to be satisfied in the sequel unless otherwise specified.

***Assumption 1***

***(1-1)*** ***Boundedness:*** *For all* $(x,y,\mu) \in D \times \Gamma \times \mathcal{M}_2^{(\Gamma)}$, *there is a constant* $\bar{U} > 0$ *such that* $0 \leq U \leq \bar{U}$.

***(1-2)*** ***Lipschitz continuity in actions:*** *For all* $(x_1,x_2,y,\mu) \in D \times D \times \Gamma \times \mathcal{M}_2^{(\Gamma)}$, *there is a constant* $L_D > 0$ *such that*

$$\left|U(x_1,y,\mu) - U(x_2,y,\mu)\right| \leq L_D |x_1 - x_2|. \tag{11}$$

***(1-3)*** ***Lipschitz continuity in action distributions 1:*** *For all* $(x,y,\mu,\nu) \in D \times \Gamma \times \mathcal{M}_2^{(\Gamma)} \times \mathcal{M}_2^{(\Gamma)}$, *there is a constant* $L_{\mathcal{M}} > 0$ *such that*

$$\left|U(x,y,\mu) - U(x,y,\nu)\right| \leq L_{\mathcal{M}} \|\mu-\nu\|_{\mathrm{TV}}^{(\Gamma)}. \tag{12}$$

***(1-4)*** ***Lipschitz continuity in action distributions 2:*** *For all* $(x,y,\mu,\nu) \in D \times \Gamma \times \mathcal{P}^{(\Gamma)} \times \mathcal{P}^{(\Gamma)}$, *there is a constant* $L_W > 0$ *such that*

$$\left|U(x,y,\mu) - U(x,y,\nu)\right| \leq L_W W_1(\mu,\nu). \tag{13}$$

***(1-5)*** ***Lipschitz continuity in types:*** *For all* $(x_1, y_1, y_2, \mu) \in D \times \Gamma \times \Gamma \times \mathcal{M}_2^{(\Gamma)}$ *, there is a constant* $L_\Gamma > 0$ *such that*

$$\left| U(x, y_1, \mu) - U(x, y_2, \mu) \right| \le L_\Gamma \left| y_1 - y_2 \right|. \tag{14}$$

***(1-6)*** ***Behavior near maxima:*** *For all* $(y, \mu) \in \Gamma \times \mathcal{M}_2^{(\Gamma)}$ *, the following holds true: the function* $U(x, y, \mu)$ *of* $x \in D$ *has a finite number of strict local maximum points* $x'$ *, and there is a constant* $c_U > 0$ *with* $U(x', \mu) - U(x, \mu) \approx c_U \left| x - x' \right|$ *when* $\left| x - x' \right|$ *is small; the difference between the maximum and minimum of* $U$ *is not smaller than a constant* $\Delta U > 0$ *.*

We need to define utilities not only for $\mathcal{P}^{(\Gamma)}$ but also for $\mathcal{M}_2^{(\Gamma)}$ due to technical reasons (generalized Picard–Lindelöf theorem); we need to consider a slightly larger space for (local) Lipschitz continuity of the right-hand side of logit dynamics (e.g., Proof of Theorem 3.4 in Lahkar and Riedel [3]), and the space $\mathcal{M}_2^{(\Gamma)}$ is a reasonable choice in our case. The assumption **(1-4)** is used only for the well-posedness of the McKean–Vlasov SDEs and is stronger than **(1-3)**. The assumptions **(1-1)** through **(1-3)** are rather standard in view of evolutionary games, including logit dynamics [1,3]. The upper and lower bounds of $U$ can be replaced by other bounded values without any loss of generality.

The assumption **(1-5)** ensures that the right-hand side, and hence solutions to logit dynamics, are Lipschitz continuous with respect to agent types. The last assumption **(1-6)**, which is the most crucial for the logit dynamic subject to exploration costs, is that the utilities are not flat (e.g., constant); indeed, dealing with a constant utility is not interesting, but the more striking fact is that this kind of variability assumption about utilities is necessary for proving the Lipschitz continuity of the solution-dependent inverse temperature (see **Section 3.1**). Without this assumption, even with the regularization term, the Lipschitz continuity seems to remain at most local, which is not enough for proving the uniqueness of solutions to our logit dynamic. Later, we numerically investigate that assumption **(1.6)** may not be necessary for stable computation of the logit dynamic.

With respect to assumptions about utilities, the difference between the previous study [26] is that we need not assume the extraordinarily strong assumption about the relative entropy of logit functions; instead, we employ a regularization method to the cost constraint. Note that a regularization method with a different form has been applied to replicator-type dynamics in Yoshioka [26], and this study shows that our regularization techniques may apply to a certain range of evolutionary game dynamics, including logit and replicator dynamics.

An example that satisfies **Assumption 1** is given as follows: for any $(x, y, \mu) \in D \times \Gamma \times \mathcal{M}_2^{(\Gamma)}$,

$$U(x, y, \mu) = \int_{x \in D} \int_{w \in \Gamma} G(x, y, w) |\mu| (\mathrm{d}x, w) \mathrm{d}w + h(x), \tag{15}$$

where $h$ is a nonnegative and nonzero function and $G : D \times \Gamma \times \Gamma \to \mathbb{R}$ is a continuous function modeling interaction among heterogeneous agents, such as a graphon [40,41]. Here, $|\mu|$ is the total variation measure of $\mu$, and defining a utility by using $|\mu|$ is not harmful in applications because $|\mu| = \mu$ if $\mu \in \mathcal{P}^{(\Gamma)}$. A more complex example would be, with some suitably smooth function $F$,

$$U(x,y,\mu) = F\left(x, y, \int_{x\in D}\int_{w\in\Gamma} G(x,y,w)|\mu|(\mathrm{d}x,w)\mathrm{d}w\right) + h(x). \tag{16}$$

### 2.3 The logit dynamic

Now, we formulate our logit dynamic. Fix $\varepsilon > 0$ and $\omega > 0$. For any $(x,y,\mu) \in D \times \Gamma \times \mathcal{M}^{(\Gamma)}$, set the logit function

$$\mathbb{L}_\beta(x,y,\mu) = \frac{\exp(\beta U(x,y,\mu))}{\int_D \exp(\beta U(z,y,\mu))\mathrm{d}z}. \tag{17}$$

For any Borel measurable set $A$ in $D$, our logit dynamic is presented as follows:

$$\frac{\mathrm{d}\mu_t(A,y)}{\mathrm{d}t} = \int_A \mathbb{L}_{\beta(y,\mu_t)}(z,y,\mu)\mathrm{d}z - \mu_t(A,y), \quad t > 0 \quad y \in \Gamma \tag{18}$$

with $\beta(y,\mu_t) > 0$ being a solution to

$$\underbrace{\int_D \left\{ \mathbb{L}_\beta(z,y,\mu_t)\ln\mathbb{L}_\beta(z,y,\mu_t) - (\mathbb{L}_\beta(z,y,\mu_t) - 1) \right\}\mathrm{d}z}_{\text{True cost}} + \underbrace{\omega\beta}_{\text{Regularization}} = \varepsilon. \tag{19}$$

There are two distinct points between the dynamics (1) and (17). First, the latter accounts for agent heterogeneity. Second, the inverse temperature is solution dependent, and this fact introduces additional nonlinearity to the problem. Yoshioka [26] studied a homogeneous logit dynamic subject to an exploration cost constraint under a more severe, entropic assumption about utilities. In contrast, this study employs **Assumption 1** only. The price to pay is the regularization term $\omega\beta$ included in (19). This term ensures that equation (19) is uniquely solvable for $\beta > 0$ and that the solution is uniformly bounded. Without this bound, we could not show the boundedness and continuity result of $\mathbb{L}_{\beta(y,\mu_t)}$ that is needed for the Picard–Lindelöf theorem.

The regularization term may be understood as a kind of barrier function [e.g., 48,49] that penalizes the extreme growth of $\beta$ although not seemingly classical, and the strength of penalization is controlled by the coefficient $\omega$. The role of this regularization term is to avoid infinitely large values of the inverse temperature (i.e., inverse of a Lagrangian multiplier), which thus acts like a barrier function method. This regularization method combined with **(1-4)** in **Assumption 1** leads to a Lipschitz continuity estimate for the right-hand side of our logit dynamic, which is the key for its well-posedness. Our regularization term $\omega\beta$ has quite a simple form, but technically other functional forms are allowed if it vanishes at $\beta = 0$ and is positive and strictly increasing for $\beta > 0$, such as $\omega(\exp(\beta) - 1)$.

### 2.4 The McKean–Vlasov SDE

Next, we present our McKean–Vlasov SDE. The low of $X_t^y$ is denoted by $\mathcal{L}\left(X_t^y\right)$. The McKean–Vlasov SDE is formulated as follows:

$$\mathrm{d}X_t^y=\int_D\int_0^\infty\left(z-X_t^y\right)\mathbb{I}_{\{w\le\mathbb{L}_{\beta(y,\mathcal{L}_{t-})}\}}N\left(\mathrm{d}t,\mathrm{d}z,\mathrm{d}w,\{y\}\right),\quad t>0 \tag{20}$$

subject to some initial condition $X_0^y\in D$ ( $y\in\Gamma$ ), where $N$ is a Poisson random measure on $[0,+\infty)\times D\times(0,+\infty)\times\Gamma$. We write $\mathcal{L}_t=\left\{\mathcal{L}\left(X_t^y\right)\right\}_{y\in\Gamma}$. We expect that there is consistency between the logit dynamic (18) and the SDE (20) because the former can be seen as an FPE of the latter; in particular, $\mathcal{L}\left(X_t^y\right)=\mu(\cdot,t)$. Here, the Poisson random measure $N$ is understood as a continuum of Poisson random measures on $[0,+\infty)\times D\times(0,+\infty)$, developed in the context of Fubini extension [e.g., 50,51].

## 3. Mathematical analysis

### 3.1 The logit dynamic

A key in our logit dynamic is the equality (19), particularly the boundedness and Lipschitz continuity with respect to agent actions of $\beta\left(y,\mu_t\right)$. The following **Lemma 1** addresses this issue. In the rest of **Section 3**, we fix $\varepsilon,\omega>0$.

***Lemma 1***

*For any $(y,\mu)\in\Gamma\times\mathcal{M}_2^{(\Gamma)}$, consider the following equation for $\beta\ge 0$:*

$$\int_D\left\{\mathbb{L}_\beta(z,y,\mu)\ln\mathbb{L}_\beta(z,y,\mu)-\left(\mathbb{L}_\beta(z,y,\mu)-1\right)\right\}\mathrm{d}z+\omega\beta=\varepsilon . \tag{21}$$

*This equation admits a unique solution $\beta=\beta(y,\mu)$ such that*

$$\left(\underline{\beta}:=\right)\frac{2\varepsilon}{\omega+\sqrt{\omega^2+2\varepsilon\bar{U}^2}}\le\beta(y,\mu)\le\left(\bar{\beta}:=\right)\frac{\varepsilon}{\omega}. \tag{22}$$

*Moreover, there is a constant $L_\beta>0$ independent of $(y,\mu,\nu)\in\Gamma\times\mathcal{M}_2^{(\Gamma)}\times\mathcal{M}_2^{(\Gamma)}$ such that*

$$\left|\beta(y,\mu)-\beta(y,\nu)\right|\le L_\beta\left\|\mu-\nu\right\|_{\mathrm{TV}}^{(\Gamma)}. \tag{23}$$

The next technical lemma shows some Lipschitz continuity about the logit function.

***Lemma 2***

*For any $(x,y,\mu,\nu)\in D\times\Gamma\times\mathcal{M}_2^{(\Gamma)}\times\mathcal{M}_2^{(\Gamma)}$, there exists a constant $C_0>0$ such that*

$$\left|\mathbb{L}_{\beta(y,\mu)}(x,y,\mu)-\mathbb{L}_{\beta(y,\nu)}(x,y,\nu)\right|\le C_0\left\|\mu-\nu\right\|_{\mathrm{TV}}^{(\Gamma)}. \tag{24}$$

On the basis essentially of **Lemma 2**, we state the unique existence result for the logit dynamic (18). The proof itself is by a direct adaptation of Proof of Theorem 3.4 in Lahkar and Riedel [3] with a modification because our model has an extra dimension about agent heterogeneity.

***Proposition 1***

*The logit dynamic (18) admits a unique solution* $(\mu_t)_{t\geq 0}$ *such that* $\mu_t \in \mathcal{P}^{(\Gamma)}$ *for all* $t \geq 0$.

***Remark 1*** For the logit dynamic (18), **(1-4)** in **Assumption 1** may be weakened as follows, which only assumes the boundedness of the total variation norm for all $y \in \Gamma$ but not its continuity. For example, if utilities have the form (15) with $G$ being strictly bounded and is piecewise continuous, then this replacement works; however, assuming that $G$ is strictly bounded is insufficient because very irregular functions (e.g., taking the value 1 at rational points and 0 otherwise) cannot be integrated against $\mu \in \mathcal{M}_2^{(\Gamma)}$. Whether the same applies to the McKean–Vlasov SDE is currently uncertain, and will be studied in future.

### 3.2 The McKean–Vlasov SDE

We present our results concerning the uniqueness and existence of strong (pathwise) solutions to the McKean–Vlasov SDE (20). Before stating these results, we present the following technical lemma to be used in **Proofs of Propositions 2-3**.

***Lemma 3***

*For any* $(x, y, \mu, \nu) \in D \times \Gamma \times \mathcal{P}^{(\Gamma)} \times \mathcal{P}^{(\Gamma)}$, *there exists a constant* $C_1 > 0$ *such that*

$$\left|\mathbb{L}_{\beta(y,\mu)}(x,\mu) - \mathbb{L}_{\beta(y,\nu)}(x,\nu)\right| \leq C_1 W_1(\mu(\cdot,y), \nu(\cdot,y)). \tag{25}$$

The following **Propositions 2-3** are a version of the uniqueness result for McKean–Vlasov SDEs driven by diffusion noise (e.g., Proof of Proposition 4 in Kuehn and Pulido [37]; Proposition 3.1 in Coppini et al. [52]); our SDE is different from existing ones because its strong solutions are bounded in $D$ because of its form and is driven by jumps. More specifically, our SDE is purely driven by jumps with finite activities, which is the reason that we can proceed with $W_1$ (see also Liang et al. [53]).

***Proposition 2***

*For each* $T > 0$, *there exists at most one strong solution to the McKean–Vlasov SDE (20) in* $\mathcal{H}_T$.

***Proposition 3***

*For each* $T > 0$, *there exists a strong solution to the McKean–Vlasov SDE (20) in* $\mathcal{H}_T$.

We conclude this section with the following **Proposition 4**, stating a consistency result between the SDE (20) and the logit dynamic (18), which is a consequence of classical Itô's formula for jump processes. This is a natural result because the McKean–Vlasov SDE (20) was motivated by the logit dynamic (18) as an FPE, but to the best of the author's knowledge there has been no explicit statement.

***Proposition 4***

*Let* $(\mu_t)_{t\geq 0}$ *be the unique solution to the logit dynamic (18) and let* $\{\mathcal{L}_t(X^{\bullet})\}_{t\geq 0}$ *be the law of the strong solution to the McKean–Vlasov SDE (20). If* $\mathcal{L}_0(X^{\bullet})=\mu_0$, *then* $\mathcal{L}_t(X^{\bullet})=\mu_t$ for all $t\geq 0$

## 4. Demonstrative applications

This section is devoted to numerical computation of the logit dynamic (18). Details of the discretization scheme, its computational procedure and convergence study are presented in **Appendix B**.

### 4.1 Application 1

We apply the logit dynamic (18) with $D=[0,1]$. Computational resolution is the computational resolution $\Delta x=\Delta y=1/400$ and $\Delta t=1/200$. This resolution is sufficiently fine for the computation here, as demonstrated in **Appendix B**. We focus particularly on the graphon $G$ in terms of the utilities and the regularization method.

The utility $U$ that we consider here has the form (16) and is explicitly specified as follows:

$$U(x,y,\mu)=\underbrace{\left(\underbrace{a_1}_{\text{Fish catch}}+\underbrace{a_2\int_\Gamma G(y,w)\int_D z\mu(\mathrm{d}z,w)}_{\text{Agent interaction}}\right)x}_{\text{Fishing utility}}-\underbrace{\frac{1}{2}x^2}_{\text{Fishing cost}}, \tag{26}$$

which is an adaptation of the utility for "Cities Game" in Carmona et al. [54] to our setting, i.e., replacing a continuum of delta functions in the interaction term with $\mu$. Here, $a_1>0$ and $a_2\in\mathbb{R}$, and $G$, a graphon, here is a continuous function on $\Gamma\times\Gamma$. We reformulated the original utility in the following two ways. First, it was originally a minimization problem of a utility, and hence, we changed its sign because our problem is a maximization problem. Second, it was originally studied by assuming agent actions as a continuum of Dirac deltas, whereas our study uses continuums of probability measures $\mu$. As indicated in (26), we provide another interpretation for the utility (26) from a fisheries perspective. In the present setting, $x$ is the (normalized) arrival intensity of agents, which are anglers, to a fishing site, and $y$ represents their types. This utility contains two terms, which are the fishing utility (first term) and the fishing cost (second term). Agent interaction is modeled in the first term through a graphon $G$ such that when an agent, an angler in the present case, is more excited about fishing if the others catch fish more when $a_2>0$, and conversely, the more successful other people are at fishing, the less enthusiastic he/she becomes when $a_2<0$. The graphon $G$ expresses the degree of agents' enthusiasm for fishing.

The linear-quadratic nature of (26) enabled its detailed mathematical analysis (Carmona et al., [54]; the unique Nash equilibrium of this utility has been formally obtained in a closed form as follows:

$$\mu(\mathrm{d}x, y) = \delta(x - \hat{x}_y), \quad (x, y) \in D \times \Gamma, \tag{27}$$

where the continuous function $\hat{x}_y$ of $y \in \Gamma$ is uniquely found from the integral equation

$$\hat{x}_y - a_2 \int_\Gamma G(y, w) \hat{x}_w \mathrm{d}w = a_1, \quad (y, w) \in \Gamma \times \Gamma, \tag{28}$$

provided that $a_2 \max_{y,w\in\Gamma} |G(y, w)| < 1$ (i.e., the right-hand side of (28) is invertible). For later use, we set the average actions:

$$\bar{x}_y = \int_D z\mu(\mathrm{d}z, y), \quad y \in \Gamma. \tag{29}$$

The formula (28) implies that the peaks of the PDFs associated with $\mu$ should be close to the curve $\hat{x}_y$ ($y \in \Gamma$), i.e., $\hat{x}_y \approx \bar{x}_y$, when $\beta$ is large or equivalently $\varepsilon$ is large. The integral equation (28) is solved by a fixed-point iteration method.

We fix $c_1 = 0.4$, $c_2 = 0.6$, and the bilinear graphon $G(y, w) = yw$, $(y, w) \in \Gamma \times \Gamma$, which is one of the most nontrivial graphons chosen here for modeling heterogeneous agent interactions. The other parameter values are set as follows unless otherwise specified: $\varepsilon = 1$ and $\omega = 10^{-4}$. These parameter values were chosen so that the heterogeneity of $\hat{x}_y$ becomes reasonably visible. Finally, we regard a numerical solution to be stationary if the PDF $p$ associated with the discretized $\mu$, which is computed by $\mu(D_{i,j}) / \Delta x$ (for the notation, see **Appendix B**), changes less than $10^{-10}$ for all the computational cells between successive time steps. The initial condition of $\mu$ is the uniform distribution in $D$ for all $y \in \Gamma$.

First, we study the transient profiles of the solutions. The computed PDFs at selected time steps are shown in **Figure 1**, which reveals their transient behavior from uniform to almost stationary with a profile that rises like a mountain range. In this case, the numerical solution becomes stationary at time 18.2. The numerical solutions obtained are smooth and nonnegative, and do not contain spurious oscillations.

Second, we investigate the equilibrium. As stated above, at a stationary state, the computed PDFs should have peaks near the Nash equilibrium $\hat{x}_y$ ($y \in \Gamma$), and this tendency should become stronger as $\beta$ decreases or $\varepsilon$ increases. The (average) equilibrium actions between the exact and computed results for $\omega = 10^{-2}$ are compared in **Figure 2**. Similarly, **Figure 3** compares the (average) equilibrium actions between the exact and computed results for $\varepsilon = 1$ and $\omega = 0$ (without any regularization). **Tables 1-2** show the maximum absolute difference between $\hat{x}_y$ and $\bar{x}_y$ for different values of $\omega$ and $\varepsilon$, respectively. According to **Figure 2**, the exact and computational results are comparable, and their maximum difference decreases as $\omega$ decreases. Choosing a sufficiently small $\omega$ results in a difference at the order of $10^{-6}$, and remarkably, the numerical computation does not breakdown even without any regularization effect, which was theoretically necessary (**Proof of Lemma 1**). This gap between theory and application is considered due to the "static" estimate employed in **Lemma 1**, where the Lipschitz continuity of $\beta$ is evaluated for any continuum of measures, whereas in each computation, we need to account for only the transient path of a numerical solution. A sharper Lipschitz estimate may be found if we evaluate $\beta$ along solutions; however, this probably necessitates a different theory from that established in this paper and will therefore be addressed elsewhere. As shown in **Table 2**, the maximum absolute difference between $\hat{x}_y$ and $\overline{x}_y$ decreases as $\varepsilon$ increases, as theoretically expected. These computational results also suggest the validity of our logit dynamic in conjunction with the numerical method because it gives a result close to the theoretical one.

We conduct an additional analysis of the balance between the two terms on the left-hand side of (19) at equilibrium:

$$\underbrace{\frac{\int_D \left\{ \mathbb{L}_{\beta(y,\mu)}(z,y,\mu) \ln \mathbb{L}_{\beta(y,\mu)}(z,y,\mu) - \left( \mathbb{L}_{\beta(y,\mu)}(z,y,\mu) - 1 \right) \right\} \mathrm{d}z}{\varepsilon}}_{\text{Cost }1(y)} + \underbrace{\frac{\omega \beta(y,\mu)}{\varepsilon}}_{\text{Cost }2(y)} = 1 . \qquad (30)$$

In this equation, $\text{Cost } 1(y)$ represents the ratio of the true exploration cost. **Table 3** shows the computed $\text{Cost } 1$, which is the average of the computed $\text{Cost } 1(y)$ for $y \in \Gamma$ with $\varepsilon = 1,2$, for different values of $\omega$ along $y = 0.5$. In this case, $\text{Cost } 1(y)$ is almost independent of $y$. According to **Tables 1 and 3**, the true cost can become significantly smaller than 1 even when there is a small error between the exact and computed equilibria, e.g., in case $\omega = 10^{-2}$. A comparison of the computed stationary PDFs for $\omega = 0$ and $\omega = 10^{-2}$ clearly reveals that they are qualitatively different, although both share a common unimodal nature, as shown in **Figure 4**. **Table 3** also shows that the ratio of the true cost decreases for a larger maximum cost $\varepsilon$. These results demonstrate that specifying an excessively large regularization parameter $\omega$ can lead to quantitatively different results for the probability density function, even if the result is reasonable for the Nash equilibrium. To avoid such situations, it would be generally preferable to use a small regularization term in numerical calculations. In addition, this problem demonstrates that computation can sometimes proceed without it.

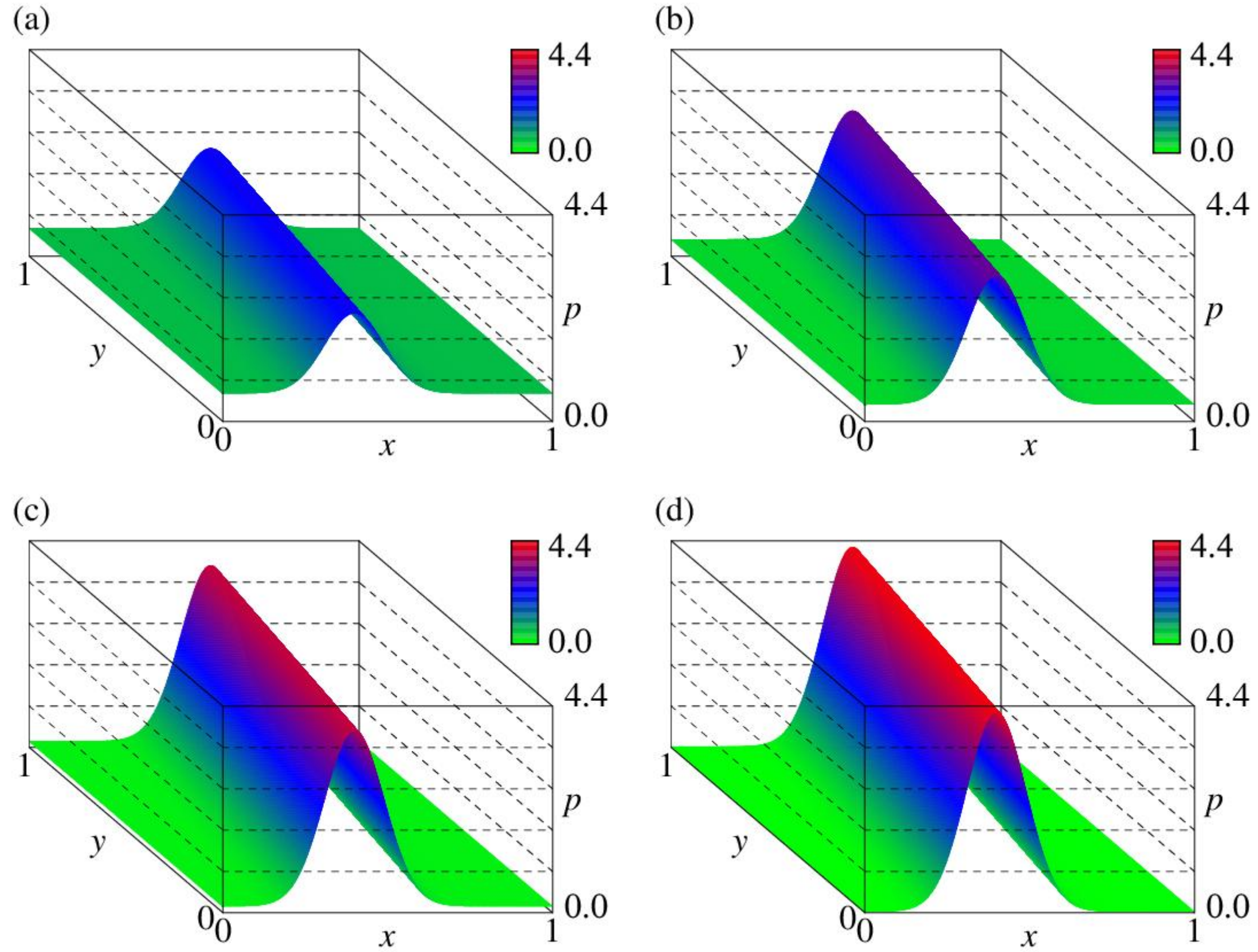

**Figure 1.** Computed PDFs at selected time steps: (a) $t = 0.5$, (b) $t = 1$, (c) $t = 2$, and (d) $t = 4$.

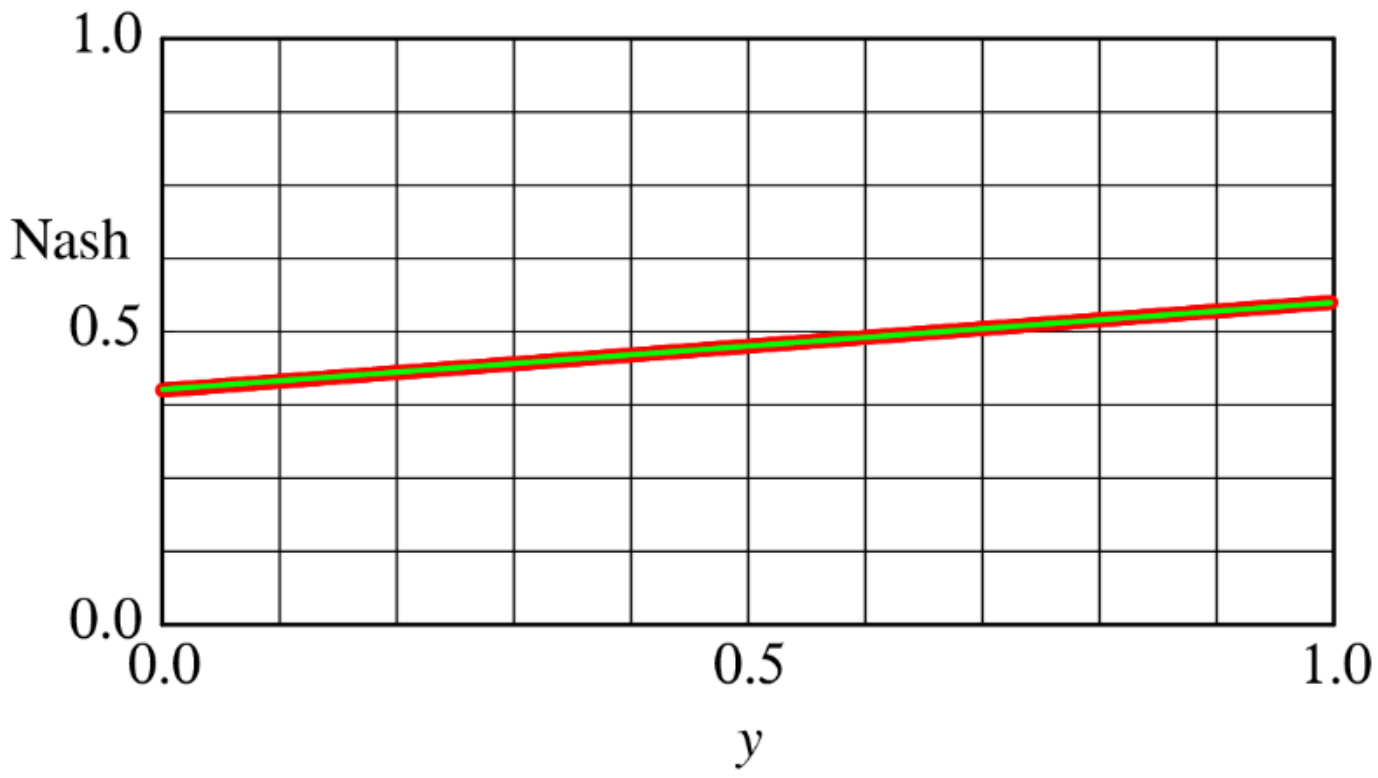

**Figure 2.** Equilibrium actions between the exact (red) and computed results with $\omega = 10^{-2}$ (green).

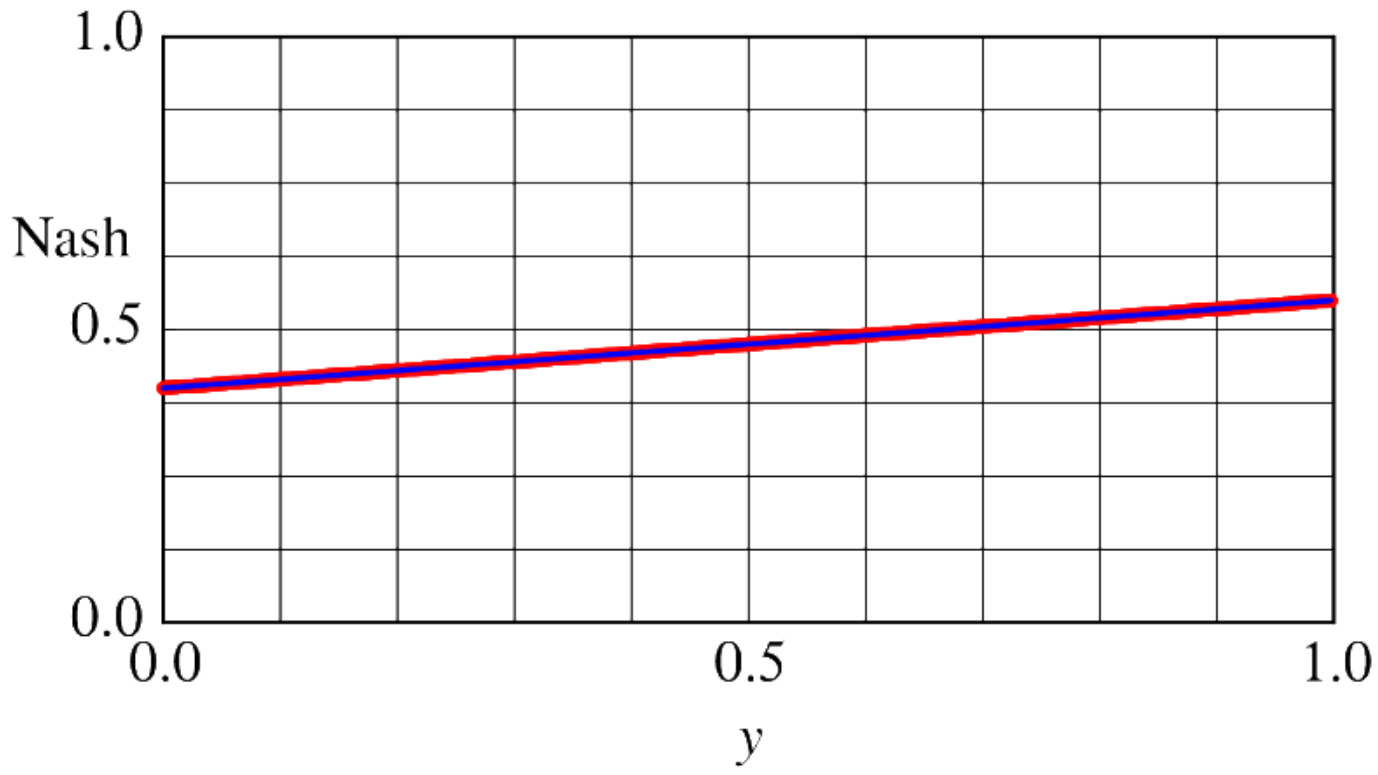

**Figure 3.** Equilibrium actions between the exact (red) and computed results with $\varepsilon = 1$ and $\omega = 0$ (blue).

**Table 1.** Maximum absolute difference between $\hat{x}_y$ and $\bar{x}_y$ for different values of $\omega$.

| $\omega$ | 0 | $10^{-4}$ | $10^{-3}$ | $10^{-2}$ |
|---|---|---|---|---|
| Difference | 1.452E-06 | 1.879E-06 | 1.054E-05 | 1.294E-03 |

**Table 2.** Maximum absolute difference between $\hat{x}_y$ and $\bar{x}_y$ for different values of $\varepsilon$, where $\omega = 0$.

| $\varepsilon$ | 0.5 | 1.0 | 1.5 | 2.0 |
|---|---|---|---|---|
| Difference | 1.652E-03 | 1.452E-06 | 3.119E-10 | 1.779E-10 |

**Table 3.** Cost 1 for different values of $\omega$.

| $\omega$ | 0 | $10^{-4}$ | $10^{-3}$ | $10^{-2}$ |
|---|---|---|---|---|
| Cost $1: \varepsilon = 1$ | 1 | 9.877E-01 | 8.973E-01 | 5.215E-01 |
| Cost $1: \varepsilon = 2$ | 1 | 9.602.E-01 | 7.948.E-01 | 4.607.E-01 |

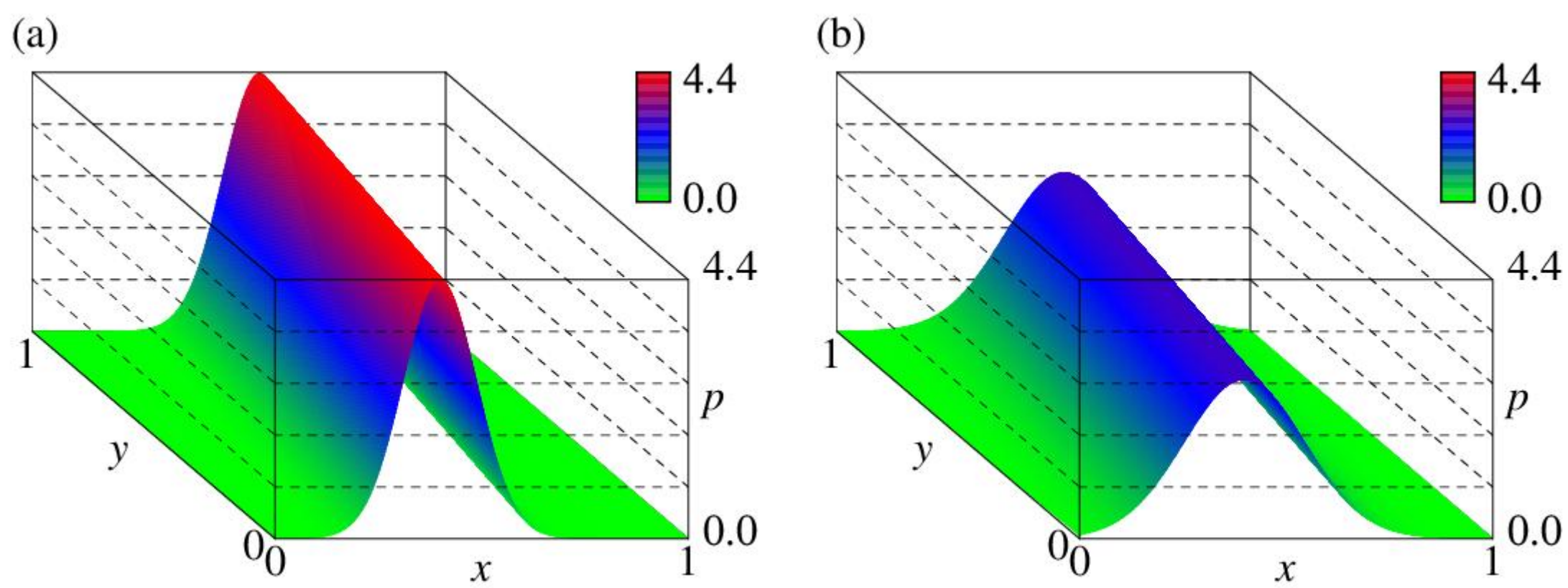

**Figure 4.** Computed stationary PDFs: (a) $\omega = 0$ and (b) $\omega = 10^{-2}$.

### 4.2 Application 2

#### 4.2.1 Problem background

We carry out another computational study of logit dynamics, which is based on the previous computational setting. The study area is the middle section of the Hii River, a first-class river that flows through eastern Shimane Prefecture, Japan, and is in the middle reach (**Figure 5**). This river reach is in a popular fishing area for *Plecoglossus altivelis* (called Ayu in Japan) in the Hii River and other fish during summer to autumn (July to October), which are major inland water resources in Japan, through angling and casting-net fishing [55]. Previous studies have shown that the spatial distribution of aquatic organisms in rivers, including Ayu, can be effectively modeled using the habitat SI [e.g., 56,57,58] based on local hydraulic variables such as water depth and flow velocity. Just upstream of this river section lies the Obara Dam, the largest multipurpose dam in the Hii River system. Currently, this dam does not have hydroelectric power generation facilities, but a plan is underway to add them by the 2030s[1]. If the water released from the dam was diverted for power generation, the amount of water flowing through this river section would decrease, which should alter the spatial distribution of the Ayu's SI, i.e., the quality of the fishing grounds. This motivated us to analyze how the presence or absence of water intake from the dam could theoretically alter anglers' behavior. Although this is a case study focusing on a single region, we believe that the proposed framework itself can be applied to other rivers and other fish species.

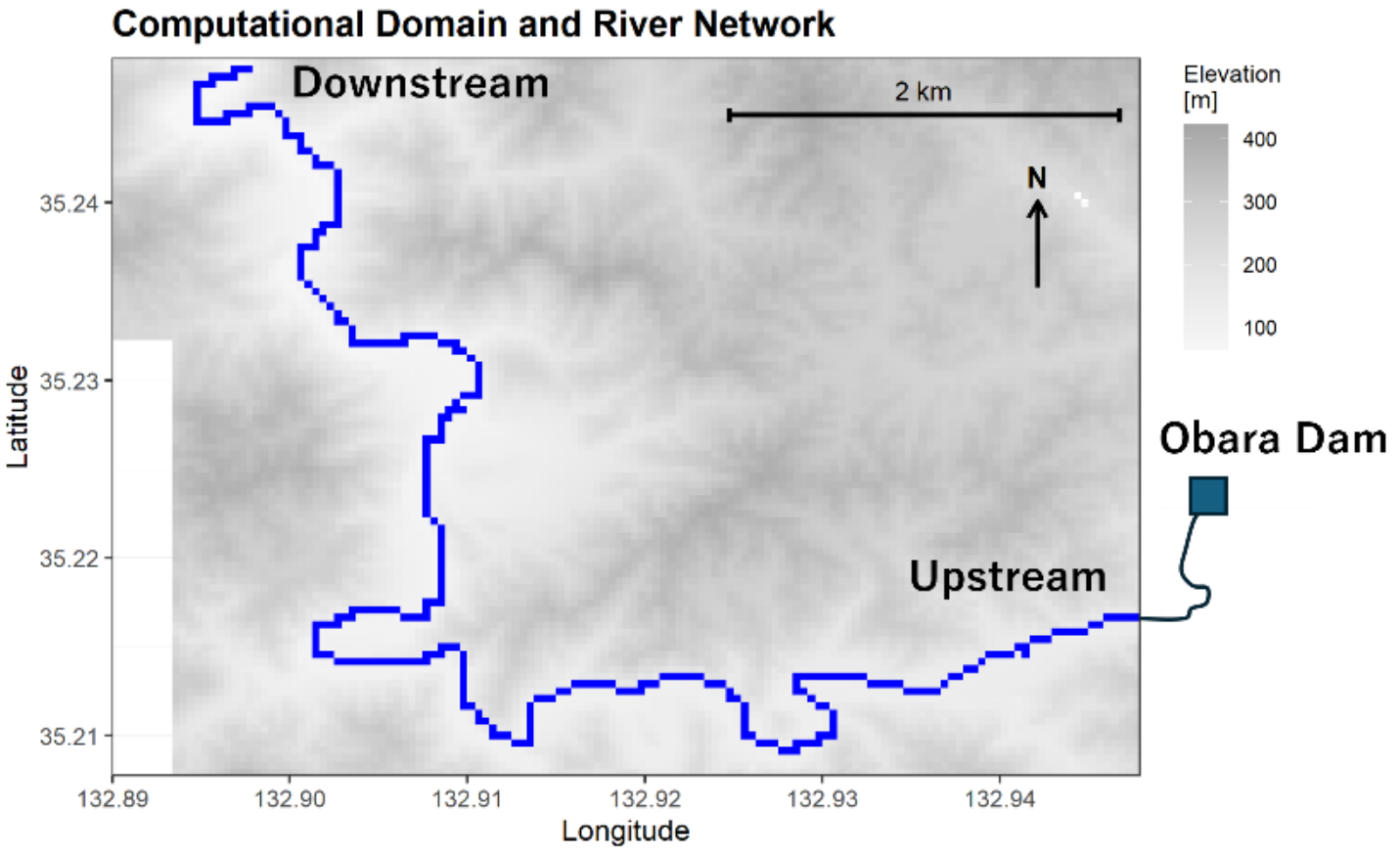

**Figure 5.** Study site. The blue line shows the one-dimensional river meshes of the Hii River in a hydrological–hydraulic model; the blue rectangle shows Obara Dam.

[1] Shimane Prefecture Press Release. https://www3.pref.shimane.jp/houdou/articles/165294 (in Japanese).

#### 4.2.2 Problem setting 1: Hydrological–hydraulic model

We use a hydrological–hydraulic model for computing SI distributions [45] that has been validated using water level observation data in the target river reach. Details regarding the model formulation and computational procedures are presented in Tanaka et al. [45]. In this model, the target river reach is treated as a one-dimensional manifold, and the water flow occurring there is simulated using spatially two-dimensional slope flow and one-dimensional local inertial equations driven by discharge from the dam. Here, local inertial equations are fundamental models for simulating surface water flows developed in hydrology [59] and have been applied to many engineering applications [60,61,62]. The SI of *P. altivelis* at each location in the river is computed from the output results of the hydrological-hydraulic model by using the methodology of Onitsuka et al. [63] following Tanaka et al. [45]: $\mathrm{SI}(V) = \alpha^{\lambda} V^{\lambda-1} e^{-\alpha V} / \Gamma(\lambda)$ with parameters $\alpha = 1.5$ and $\lambda = 1.25$, where $V \geq 0$ is the flow speed in (m/s) and $\Gamma(\cdot)$ is the gamma function. This empirical SI depends only on flow velocity and is applicable only to *P. altivelis* but is easy to implement. Note that the problem presented here is not for precise predictions of the angler distribution but rather to discuss how our logit dynamic can be used in applications.

In the computation, we use the hourly output data from the hydrological–hydraulic model. The total number of computational cells corresponding to the target river reach is 302, and the computational period is from the beginning of 2019 to the end of this year, during which the model has been calibrated against the observation data of Tanaka et al. [45]. The three cases below consider possible water intake for investigating the influences of river hydrology on logit dynamics: the baseline case without water intake, water intake with an amount of 5 ($m^3$/s) from the dam (called the 5-reduced case), and water intake with an amount of 10 ($m^3$/s) from the dam (called the 10-reduced case). In the second and third cases, the amount of water intake is truncated so that the dam discharge does not fall below 1 ($m^3$/s), which is the minimum discharge from the dam. The dam discharge, water velocity, and SI for each case are shown in **Figure 6**. Discharges are significantly larger than the thresholds values during flood events, and the 10-reduced case is the extreme case where the water is almost always abstracted for the hydropower generation.

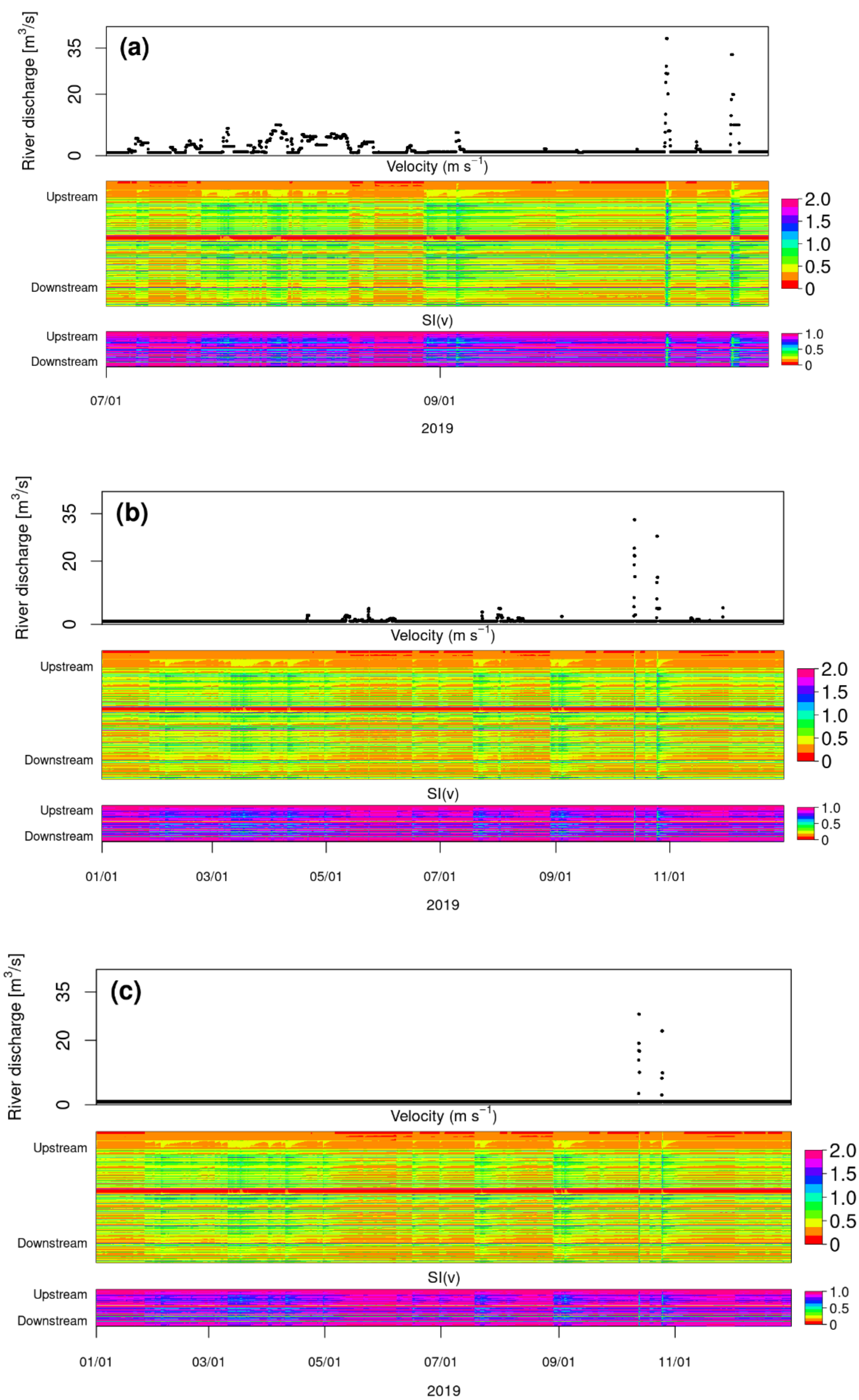

**Figure 6**: The dam discharge, water velocity, and SI for each case: (a) baseline case, (b) 5-reduced case, and (c) 10-reduced case.

#### 4.2.3 Problem setting 2: Utility and dynamics

The river reach is normalized to the closed interval $\Gamma = [0,1]$, where $y = 0$ and $y = 1$ correspond to the upstream-end and downstream-end, respectively. The other closed interval $D = [0,1]$ is considered the arrival intensity of anglers to the river in the unit (1/day) during the main fishing season (July 1 to October 31). Here, we assume that each angler arrives at the river at most once a day. Moreover, because the SI depends on time $t$, we allow for the time-dependence of graphons on $t$, which does not affect the well-posedness result in this paper if the dependence is piecewise continuous in $t$ and $y$, which is assumed in the sequel (see also **Remark 1**). Indeed, utilities considered below is piecewise constant in the $y$-direction because the flow velocity available in the hydrological-hydraulic model is a cell-wise (i.e., piecewise-constant) variable. We also assume that the reactions of anglers to river environmental conditions occur on daily basis. The unit time in the logit dynamic is then set to 1 (day). The time increment is set as $\Delta t = 1/24$ (day). The resolution is $\Delta x = \Delta y = 1/400$. The initial condition is against the uniform distribution specified at time 0 (the beginning of January 1), but the influence of the initial condition disappears before the opening time of fishing (approximately June to July, which is more than 150 days).

The utility examined here is the following, which is a version of (26) accounting for the SI $S_t(y)$ at the location $y$ at time $t$ evaluated based on the flow velocity from the nearest computational cell of the hydrological-hydraulic model:

$$U_t(x,y,\mu) = \underbrace{\underbrace{S_t(y)}_{\text{Habitat suitability}} \left( \underbrace{a_1}_{\text{Fish catch}} + \underbrace{a_2 \int_\Gamma G(y,w) \int_D z\mu(\mathrm{d}z,w)}_{\text{Agent interaction}} \right) x}_{\text{Fishing utility}} - \underbrace{\frac{1}{2}x^2}_{\text{Fishing cost}}, \quad t > 0, \tag{31}$$

where $x$ is the arrival intensity (1/day) of anglers. Here, for simplicity, we assume that the SI proportionally affects the fishing utility in a way that anglers who go to areas with high SI values potentially obtain higher utilities. A more complex dependence can be considered theoretically and addressed computationally, but we do not have such data; thus, we examine the simplest case here. The other parameter values in the utility are the same as those in the previous application example. Finally, within our framework, $\varepsilon$ can be understood as the exploration cost of finding a good fishing spot in the target river reach. In this context, the cost based on the entropy in the logit dynamic can be linked to the uncertainty in the information regarding the SI of Ayu.

#### 4.2.4 Results and discussion

We consider the following three graphons: constant graphon $G(y,w) \equiv 1$, bilinear graphon $G(y,w) = yw$, and the SI graphon $G(y,w) = S_t(y)S_t(w)$. The uniform graphon is the simplest case where agents are homogeneous. The bilinear graphon assumes that the Ayu fishing is more exciting in more downstream areas for anglers, which is an assumption made for simplicity, but a possible explanation is that a more downstream reach is closer to city areas. The SI graphon case assumes that the Ayu fishing is more exciting in areas with higher SI values, where the existence of the fish is more likely. The identification of graphon itself is an interesting topic related to social sciences, but it remains a challenge for the future. However, based on the mathematical and numerical methods presented in this paper, we believe that even complex graphons can be efficiently dealt with.

We set the exploration cost $\varepsilon = 1.5$ unless otherwise specified and examine other values of $\varepsilon$ in this subsection and investigate smaller values of $\varepsilon$ later. Case $\varepsilon = 1.5$ yields $\beta = O\left(10^2\right)$, which is significantly greater than 1. We set the regularization parameter value $\omega = 10^{-10}$; reducing this value does not affect the computational results below, and furthermore, the average Newton iteration at each time step and each $y$ grid point is approximately 3 to 5 for each computational case, suggesting the rapid convergence of the numerical algorithm at each time step. The computational time is around 10 minutes by using a common laptop for the full computation from the beginning of January 1 to the end of December 31.

**Figure 7** shows the computed PDFs of the arrival intensity at the beginnings on July 1, August 1, September 1, and October 1 with the SI graphon for the baseline case. Similarly, **Figures 8-9** show the results for the 5-reduced and 10-reduced cases. The irregularity of the PDFs in the $y$-direction is due not to computational instabilities but to the corresponding spatial distribution of the SI, as shown in **Figure 6**. Another reason for this phenomenon is that the proposed logit dynamic has no smoothing mechanism, e.g., diffusion terms, in the $y$-direction. Moreover, the numerical solutions are considered sufficiently converged at the present resolution because doubling the resolution in both space and time does not affect the numerical solution obtained (see also **Figure B1**). Although not presented here, the PDFs for the 5- and 10-reduced cases have qualitatively the same irregular profiles.

**Figure 10** compares the average arrival intensity (AAI) over all the agent types for the baseline case, 5-reduced case, and 10-reduced case. The common finding to all cases is that increasing the water intake from the dam increases the AAI of anglers, irrespective of the graphons. Agent actions, even on average, fluctuate depending on the SI dynamics during the fishing season. This is considered due to the decreased discharge from the dam, which increased the area where the SI value of Ayu through the water velocity, and this is by no means a general assertion that water intake from the dam for power generation improves the habitat of Ayu. The analysis results suggest that the assumption that anthropogenic environmental modifications always have a negative impact on the ecology of target species is not necessarily true.

We compare the influences of the exploration cost $\varepsilon$ on the average actions, by focusing on the 5-reduced case with the SI graphon. **Figure 11** shows the AAIs for different values of $\varepsilon$. As the value of $\varepsilon$ decreases, i.e., as the upper limit of the exploration cost decreases, AAI tends to approach 0.5, which corresponds to the equilibrium AAI without any heterogeneity under the uniform distribution in $x$. Indeed, the average absolute deviations between the AAI and 0.5 are as follows: 0.0459 ($\varepsilon = 1.5$), 0.0422 ($\varepsilon = 0.5$), 0.0222 ($\varepsilon = 0.1$), and 0.0069 ($\varepsilon = 0.01$). In other words, from the perspective of this study, the act of incurring exploration costs to find good fishing grounds can be the root cause of the spatiotemporal fluctuations in the agent actions.

Finally, a more complex utility can be considered, where angers avoid fishing at locations where the flow velocity is high. In this case, we can drop the first term in (31) when the velocity at $y$ equals or exceeds some threshold value, e.g., 1.0 (m/s), which does not lose the well-posedness of the logit dynamic because the flow velocity is piece-wise constant variable as for the SI. For the baseline case, **Figure 12(a)** shows the computed PDF at the beginning of July 1 with the threshold value of 1.0 (m/s) and **Figure 12(b)** compares the AAIs with and without the utility truncation based on flow velocity, and the latter with the truncation threshold of 1.0 (m/s) and 0.5 (m/s). The PDF is concentrated around where the flow velocity is small (see also **Figure 6(b)**). The AAIs are reduced by the utility truncation over the computational period, while its profile, the average angler behavior, does not qualitatively change. In the present example, the AAI reduction is more significant for the truncation threshold of 0.5 (m/s).

As demonstrated in this paper, the mathematical model and numerical methods presented in this paper make it possible to examine the agent dynamics under a variety of conditions. Although not included in this study, it is unknown whether the SI of aquatic organisms other than ayu increases due to water intake. Although our study is still in its nascent stages, it suggests the possibility of examining the relationship between changes in the river environment and the management of fishery resources from the perspective of evolutionary games.

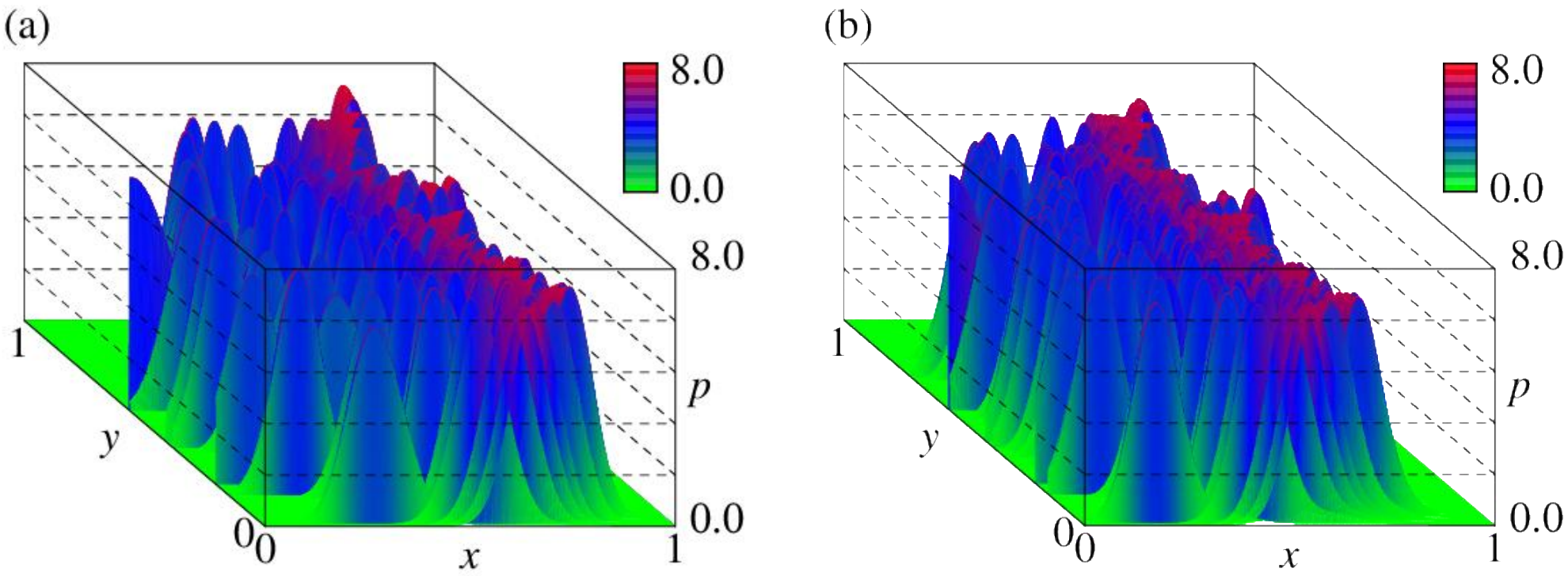

**Figure 7.** Computed PDFs of the arrival intensity at the beginnings on (a) July 1 and (b) September 1: with a uniform graph in the baseline case.

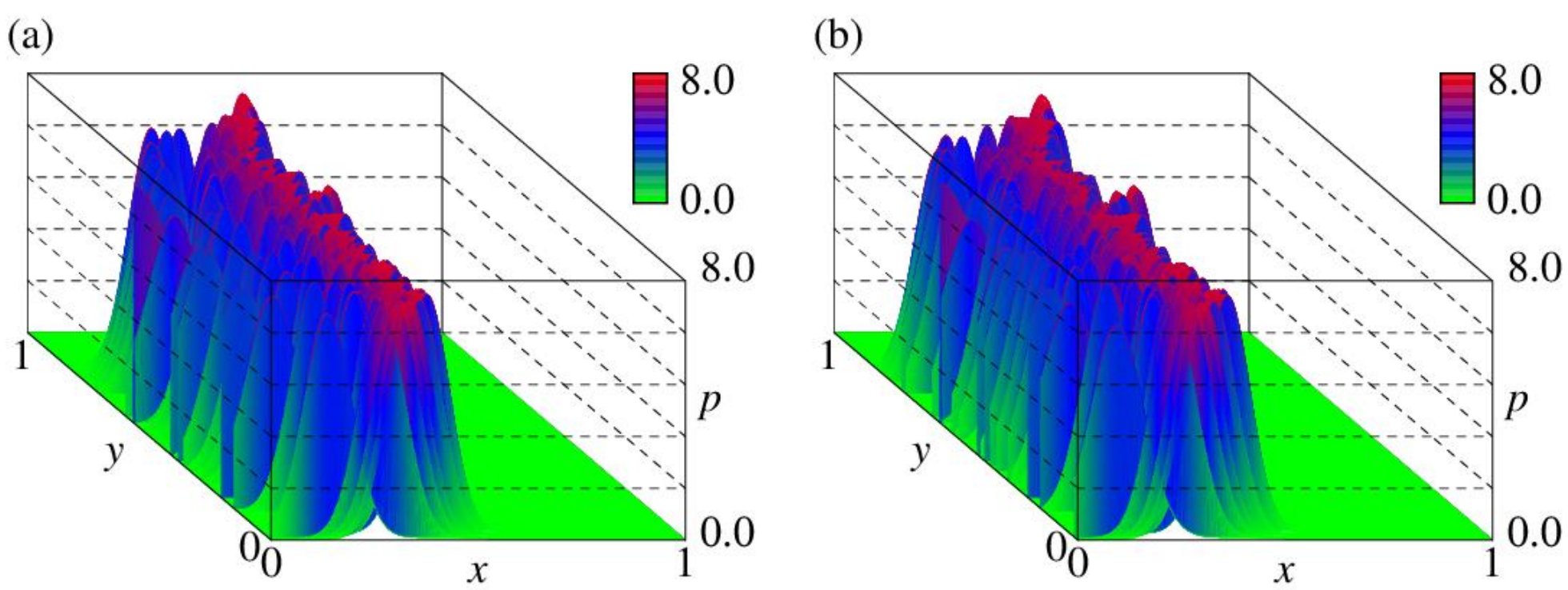

**Figure 8.** Same as **Figure 7** but with a bilinear graphon in the baseline case.

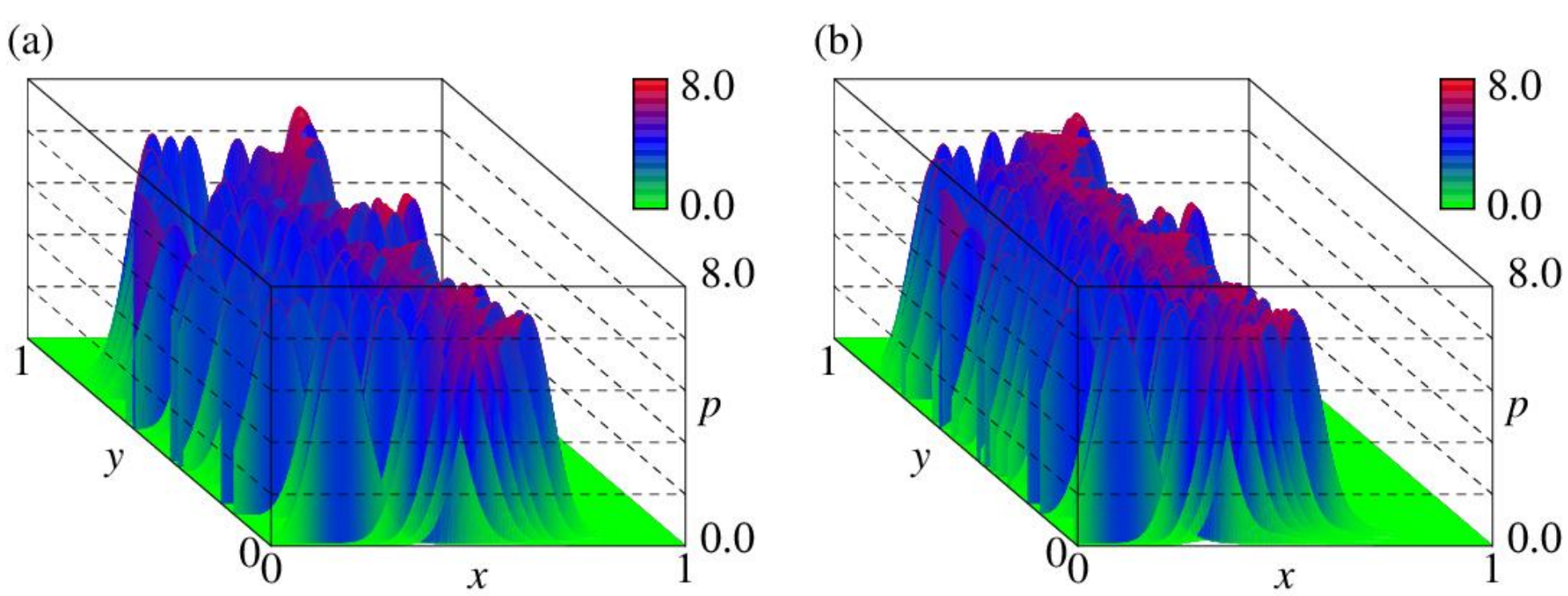

**Figure 9.** Same as **Figure 7** but with the SI graphon in the baseline case.

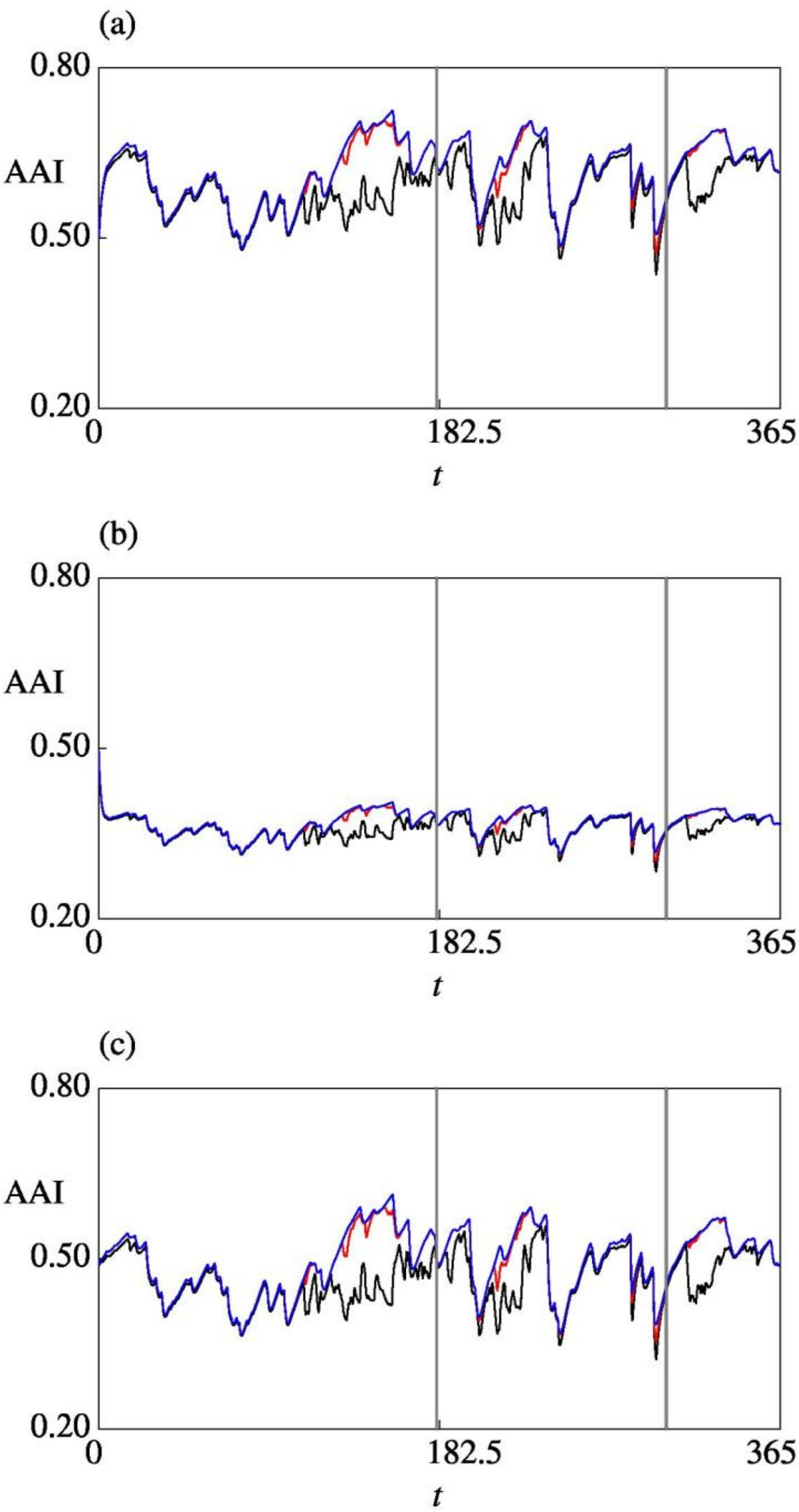

**Figure 10** Computed AAIs (1/day) in one year $0 \le t \le 365$ (day) over all the agent types for the baseline case (black), 5-reduced case (red), and 10-reduced case (blue): (a) uniform graphon, (b) bilinear graphon, and (c) SI graphon. The gray vertical lines correspond to the beginnings of July 1st and November 1st.

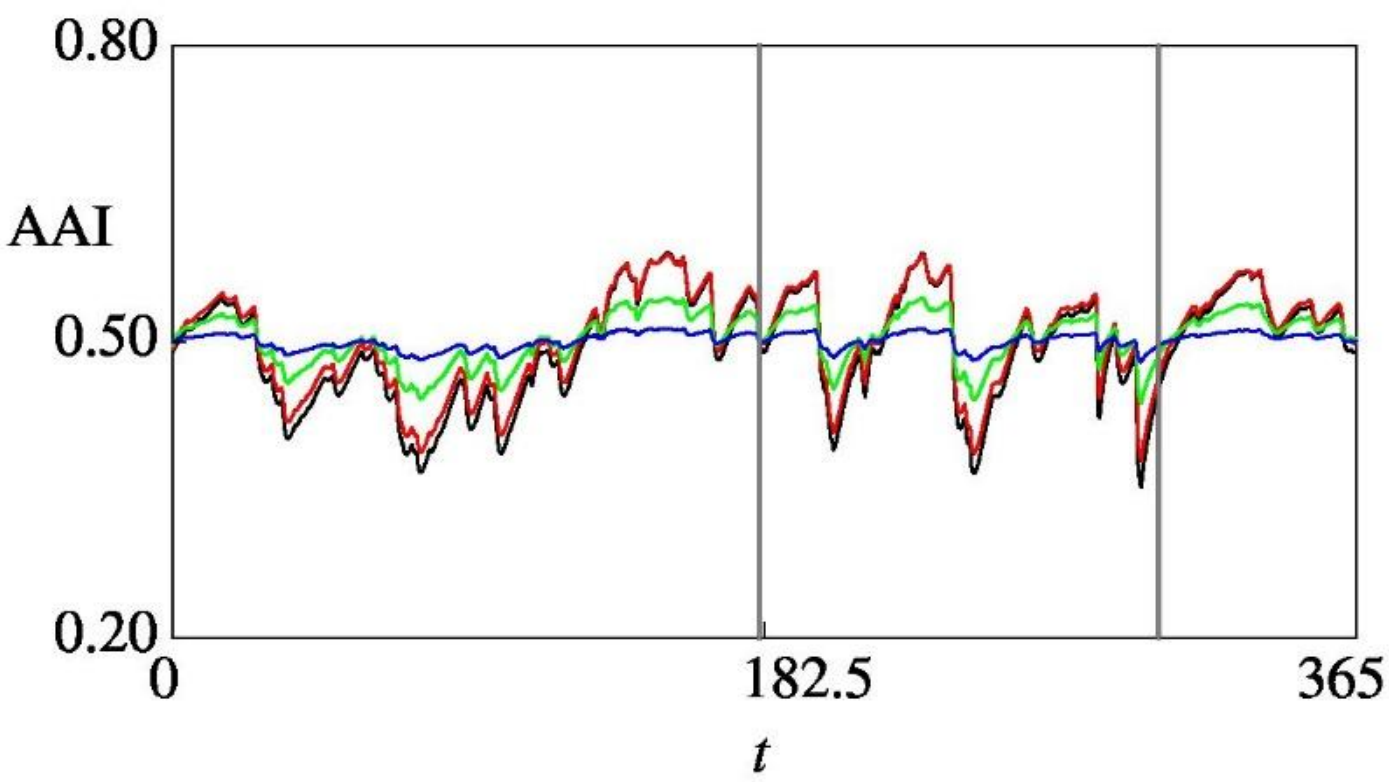

**Figure 11** Computed AAIs (1/day) in one year $0 \leq t \leq 365$ (day) over all the agent types for the SI graphon with different values of $\varepsilon$: $\varepsilon = 1.5$ (black), $\varepsilon = 1$ (red), $\varepsilon = 0.5$ (green), and $\varepsilon = 0.25$ (blue). The gray vertical lines correspond to the beginnings of July 1st and November 1st.

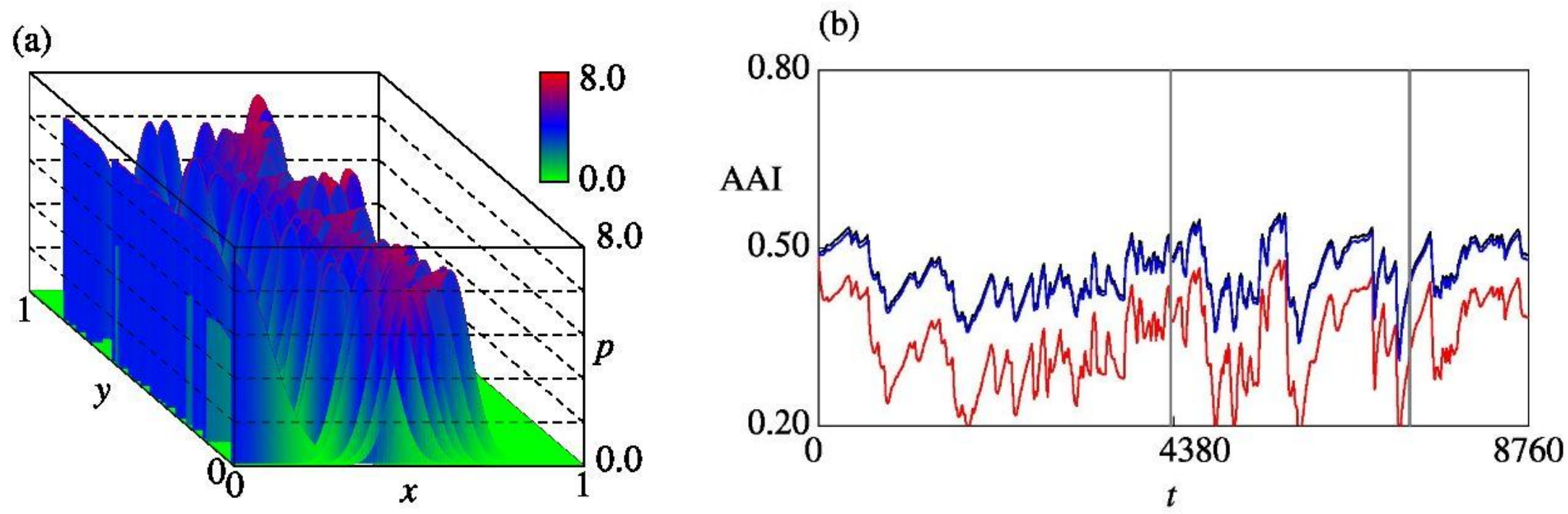

**Figure 12.** (a) Computed PDF at the beginning of July 1 with the threshold of 1.0 (m/s) (b) AAIs without velocity truncation (black), with the truncation threshold of 1.0 (m/s) (blue), and with the truncation threshold of 0.5 (m/s) (red).

## 5. Conclusion

This paper focused on logit dynamics when agents are heterogeneous and performed mathematical analyses of the FPE and associated Macken-Vlasov SDE. We showed that the unique existence of the solutions to both equations follows when a regularization term is introduced into the exploration cost constraint while imposing a specific Lipschitz continuity on utility. We computationally investigated the FPE and demonstrated that the regularization term is not innocuous in practice, but imposing a too much regularization potentially affects the resulting dynamics. The demonstrative application to a fisheries problem in a mountainous river in Japan, although under a simplified setting, suggested the broad applicability of evolutionary game theory based on nonexchangeable agents. For example, one can consider a fisheries problem in a stream network [e.g., 64].

There are several limitations in this paper. For example, we did not discuss the propagation of chaos for investigating the Macken-Vlasov SDE from a finite-agent model. This issue can be critical when an efficient and convergent discretization scheme for the SDE is constructed. Also, we did not address mean-field game versions of the logit dynamic, which potentially become a more versatile model to represent a diverse social phenomenon. In addition, our contribution is limited in the sense that we investigated only logit dynamics, while many other evolutionary game models exist in the literature. Their mathematical properties and behavior, particularly those of the Macken-Vlasov SDE, have not been extensively studied, and much room remains for further research. The development of theories with more irregular graphons is also important. Numerical computations of Macken-Vlasov SDE become important when viewing logit dynamics as consensus models based on nonconcave utilities. For high-dimensional problems, an SDE-based approach may be more efficient. With respect to practical applications of evolutionary game models when agents exhibit heterogeneity, there is still a lack of accumulated research knowledge. Even a simple application example can be considered a step toward mastering the practical applications of mathematical models. Theoretical convergence study of the proposed numerical method also remains unresolved. Currently, the authors are conducting applied studies on energy and food resource management based on evolutionary game theory.

## Appendix

### A. Proofs

***Proof of Lemma 1***

We fix $(y,\mu)\in\Gamma\times\mathcal{M}_2^{(\Gamma)}$. The first term on the left-hand side of (21) is denoted by $F(\beta)$. We have $F(0)=0$. We also have $\lim_{\beta\to+\infty}\left(F(\beta)+\omega\beta\right)=+\infty$. The upper bound in (22) then follows from the nonnegativity of $F(\beta)$. In the rest of the proof, the expectation with respect to $\mathbb{L}_\beta$ is denoted by $\mathbb{E}_\beta$ to reduce notations.

Next, elementary calculations show

$$\begin{aligned}\frac{\mathrm{d}F(\beta)}{\mathrm{d}\beta}&=\frac{\mathrm{d}}{\mathrm{d}\beta}\int_D\left\{\mathbb{L}_\beta(z,y,\mu)\ln\mathbb{L}_\beta(z,y,\mu)-\left(\mathbb{L}_\beta(z,y,\mu)-1\right)\right\}\mathrm{d}z\\&=\int_D\frac{\mathrm{d}\mathbb{L}_\beta(z,y,\mu)}{\mathrm{d}\beta}\ln\mathbb{L}_\beta(z,y,\mu)\mathrm{d}z\\&=\int_D\frac{\mathrm{d}\mathbb{L}_\beta(z,y,\mu)}{\mathrm{d}\beta}\left(\beta U(z,y,\mu)-\ln\int_D e^{\beta U(m,y,\mu)}\mathrm{d}m\right)\mathrm{d}z,\end{aligned}\tag{32}$$

where the exchange between integration and differentiation is allowed due to the positivity, smoothness, and boundedness of $\mathbb{L}_\beta$ and its derivative. Another elementary calculation shows

$$\frac{\mathrm{d}\mathbb{L}_\beta(z,y,\mu)}{\mathrm{d}\beta}=\frac{U(z,y,\mu)e^{\beta U(z,y,\mu)}}{\int_D e^{\beta U(m,y,\mu)}\mathrm{d}m}-\frac{e^{\beta U(z,y,\mu)}\int_D U(m,y,\mu)e^{\beta U(m,y,\mu)}\mathrm{d}m}{\left(\int_D e^{\beta U(m,y,\mu)}\mathrm{d}m\right)^2}.\tag{33}$$

Afterward, we obtain

$$\begin{aligned}&\int_D\beta U(z,y,\mu)\left(\frac{U(z,y,\mu)e^{\beta U(z,y,\mu)}}{\int_D e^{\beta U(m,y,\mu)}\mathrm{d}m}-\frac{e^{\beta U(z,y,\mu)}\int_D U(m,y,\mu)e^{\beta U(m,y,\mu)}\mathrm{d}m}{\left(\int_D e^{\beta U(m,y,\mu)}\mathrm{d}m\right)^2}\right)\mathrm{d}z\\&=\beta\int_D\frac{\left(U(z,y,\mu)\right)^2e^{\beta U(z,y,\mu)}}{\int_D e^{\beta U(m,y,\mu)}\mathrm{d}m}\mathrm{d}z-\beta\left(\int_D\left(\frac{U(z,y,\mu)e^{\beta U(z,y,\mu)}}{\int_D e^{\beta U(m,y,\mu)}\mathrm{d}m}\right)\mathrm{d}z\right)^2\\&=\beta\mathbb{V}_\beta\left[U(z,y,\mu)\right]\end{aligned}\tag{34}$$

and

$$\int_D\left(-\ln\int_D e^{\beta U(m,y,\mu)}\mathrm{d}m\right)\left(\frac{U(z,y,\mu)e^{\beta U(z,y,\mu)}}{\int_D e^{\beta U(m,y,\mu)}\mathrm{d}m}-\frac{e^{\beta U(z,y,\mu)}\int_D U(m,y,\mu)e^{\beta U(m,y,\mu)}\mathrm{d}m}{\left(\int_D e^{\beta U(m,y,\mu)}\mathrm{d}m\right)^2}\right)\mathrm{d}z=0,\tag{35}$$

where $\mathbb{V}_\beta\left[U(z,y,\mu)\right]$ represents the variance of $U(z,y,\mu)$ seen as a function of a random variable $z$ with respect to the PDF $\mathbb{L}_\beta$. Consequently, we obtain

$$\frac{\mathrm{d}F(\beta)}{\mathrm{d}\beta}=\beta\mathbb{V}_\beta\left[U(z,y,\mu)\right]\geq 0,\tag{36}$$

showing that $F(\beta)$ is not decreasing. The classical intermediate value theorem along with $\omega > 0$ asserts that the equation (21) admits a unique solution $\beta(y,\mu)$.

We next prove (22) and (23). We fix $(y,\mu,\nu) \in \Gamma \times \mathcal{M}_2^{(\Gamma)} \times \mathcal{M}_2^{(\Gamma)}$. The upper bound in (22) is due to the nonnegativity of $F$. We explicitly write the dependence of $F$ on $\mu$ as $F_\mu$. We set $\hat{F}_\mu(\beta) = F_\mu(\beta) + \omega\beta$. We have

$$\left|\hat{F}_\mu(\beta(y,\mu)) - \hat{F}_\mu(\beta(y,\nu))\right| = \frac{\mathrm{d}\hat{F}_\mu(\hat{\beta})}{\mathrm{d}\beta}\left|\beta(y,\mu) - \beta(y,\nu)\right| = \left(\frac{\mathrm{d}F_\mu(\hat{\beta})}{\mathrm{d}\beta} + \omega\right)\left|\beta(y,\mu) - \beta(y,\nu)\right|, \quad (37)$$

where $\hat{\beta}$ is between $\beta(y,\mu)$ and $\beta(y,\nu)$. This implies

$$\left|\beta(y,\mu) - \beta(y,\nu)\right| = \left(\frac{\mathrm{d}F_\mu(\hat{\beta})}{\mathrm{d}\beta} + \omega\right)^{-1}\left|\hat{F}_\mu(\beta(y,\mu)) - \hat{F}_\mu(\beta(y,\nu))\right|. \quad (38)$$

Now, (36) shows that

$$F_\mu(\beta) = \int_0^\beta \gamma \mathbb{V}_\gamma\left[U(z,y,\mu)\right]\mathrm{d}\gamma \quad (39)$$

and hence

$$\varepsilon = \int_0^{\beta(y,\mu)} \gamma \mathbb{V}_\gamma\left[U(z,y,\mu)\right]\mathrm{d}\gamma + \omega\beta(y,\mu), \quad (40)$$

from which we obtain

$$\varepsilon \le \bar{U}^2 \int_0^{\beta(y,\mu)} \gamma \mathrm{d}\gamma + \omega\beta(y,\mu) = \frac{1}{2}\bar{U}^2\left(\beta(y,\mu)\right)^2 + \omega\beta(y,\mu), \quad (41)$$

showing that there is a lower bound of $\beta(y,\mu)$ independent of $\mu$; there is a constant $\underline{\beta} > 0$ such that $\beta(y,\mu) \ge \underline{\beta} = 2\varepsilon / \left(\omega + \sqrt{\omega^2 + 2\varepsilon\bar{U}^2}\right)$. As a byproduct, we obtain the lower bound in (22).

We evaluate each term in (38). We first evaluate the minimum of $\dfrac{\mathrm{d}F_\mu(\hat{\beta})}{\mathrm{d}\beta}$:

$$\frac{\mathrm{d}F_\mu(\hat{\beta})}{\mathrm{d}\beta} = \hat{\beta}\mathbb{V}_{\hat{\beta}}\left[U(z,y,\mu)\right] \ge \underline{\beta}\mathbb{V}_{\hat{\beta}}\left[U(z,y,\mu)\right]. \quad (42)$$

Therefore, we obtain

$$\left|\beta(y,\mu) - \beta(y,\nu)\right| \le \left(\underline{\beta}\mathbb{V}_{\hat{\beta}}\left[U(z,y,\mu)\right] + \omega\right)^{-1}\left|\hat{F}_\mu(\beta(y,\mu)) - \hat{F}_\mu(\beta(y,\nu))\right|. \quad (43)$$

Second, we have

$$\begin{aligned}
&\left|\hat{F}_\mu(\beta(y,\mu)) - \hat{F}_\mu(\beta(y,\nu))\right| \\
&= \left|\hat{F}_\mu(\beta(y,\mu)) - \hat{F}_\nu(\beta(y,\nu)) + \hat{F}_\nu(\beta(y,\nu)) - \hat{F}_\mu(\beta(y,\nu))\right| \\
&\le \omega\left|\beta(y,\mu) - \beta(y,\nu)\right| + \left|F_\nu(\beta(y,\nu)) - F_\mu(\beta(y,\nu))\right|.
\end{aligned} \quad (44)$$

We also have

$$
\begin{aligned}
&\left|F_\nu\left(\beta(y,\nu)\right)-F_\mu\left(\beta(y,\nu)\right)\right| \\
&\le \int_D \left|\mathbb{L}_{\beta(y,\nu)}(z,y,\nu)\ln\mathbb{L}_{\beta(y,\nu)}(z,y,\nu)-\mathbb{L}_{\beta(y,\nu)}(z,y,\mu)\ln\mathbb{L}_{\beta(y,\nu)}(z,y,\mu)\right|\mathrm{d}z \\
&\le \sup_{x\in D,y\in\Gamma}\left(\ln\mathbb{L}_{\beta(y,\nu)}(x,y,\nu)\vee\ln\mathbb{L}_{\beta(y,\nu)}(x,y,\mu)+1\right)\times\int_D\left|\mathbb{L}_{\beta(y,\nu)}(z,y,\nu)-\mathbb{L}_{\beta(y,\nu)}(z,y,\mu)\right|\mathrm{d}z \\
&\le \left(\ln\left(\exp\left(\frac{\varepsilon\bar{U}}{\omega}\right)\right)+1\right)\int_D\left|\mathbb{L}_{\beta(y,\nu)}(z,y,\nu)-\mathbb{L}_{\beta(y,\nu)}(z,y,\mu)\right|\mathrm{d}z \\
&=\left(\frac{\varepsilon\bar{U}}{\omega}+1\right)\int_D\left|\mathbb{L}_{\beta(y,\nu)}(z,y,\nu)-\mathbb{L}_{\beta(y,\nu)}(z,y,\mu)\right|\mathrm{d}z
\end{aligned}
\tag{45}
$$

and

$$
\begin{aligned}
\int_D\left|\mathbb{L}_{\beta(y,\nu)}(z,y,\nu)-\mathbb{L}_{\beta(y,\nu)}(z,y,\mu)\right|\mathrm{d}z &= \int_D\left|\frac{e^{\beta(y,\nu)U(z,y,\nu)}}{\int_D e^{\beta(y,\nu)U(m,y,\nu)}\mathrm{d}m}-\frac{e^{\beta(y,\nu)U(z,y,\mu)}}{\int_D e^{\beta(y,\nu)U(m,y,\mu)}\mathrm{d}m}\right|\mathrm{d}z \\
&\le \int_D\left|\frac{e^{\beta(y,\nu)U(z,y,\nu)}}{\int_D e^{\beta(y,\nu)U(m,y,\nu)}\mathrm{d}m}-\frac{e^{\beta(y,\nu)U(z,y,\nu)}}{\int_D e^{\beta(y,\nu)U(m,y,\mu)}\mathrm{d}m}\right|\mathrm{d}z \\
&\quad+\int_D\left|\frac{e^{\beta(y,\nu)U(z,y,\nu)}}{\int_D e^{\beta(y,\nu)U(m,y,\mu)}\mathrm{d}m}-\frac{e^{\beta(y,\nu)U(z,y,\mu)}}{\int_D e^{\beta(y,\nu)U(m,y,\mu)}\mathrm{d}m}\right|\mathrm{d}z \\
&\le 2e^{\bar{\beta}\bar{U}}\int_D\left|e^{\beta(y,\nu)U(z,y,\nu)}-e^{\beta(y,\nu)U(z,y,\mu)}\right|\mathrm{d}z \\
&\le 2\bar{\beta}e^{2\bar{\beta}\bar{U}}\int_D\left|U(z,y,\nu)-U(z,y,\mu)\right|\mathrm{d}z \\
&\le 2\bar{\beta}e^{2\bar{\beta}\bar{U}}L_{\mathcal{M}}\left\|\mu-\nu\right\|_{\mathrm{TV}}^{(\Gamma)}.
\end{aligned}
\tag{46}
$$

Consequently, we obtain

$$
\begin{aligned}
&\left|\beta(y,\mu)-\beta(y,\nu)\right| \\
&\le\left(\underline{\beta}\mathbb{V}_{\hat{\beta}}\left[U(z,y,\mu)\right]+\omega\right)^{-1}\left(\omega\left|\beta(y,\mu)-\beta(y,\nu)\right|+\left(\frac{\varepsilon\bar{U}}{\omega}+1\right)2\bar{\beta}e^{2\bar{\beta}\bar{U}}L_{\mathcal{M}}\left\|\mu-\nu\right\|_{\mathrm{TV}}^{(\Gamma)}\right)
\end{aligned}
\tag{47}
$$

and hence

$$
\begin{aligned}
\left|\beta(y,\mu)-\beta(y,\nu)\right| &\le \frac{1}{\underline{\beta}\mathbb{V}_{\hat{\beta}}\left[U(z,y,\mu)\right]}\left(\frac{\varepsilon\bar{U}}{\omega}+1\right)2\bar{\beta}e^{2\bar{\beta}\bar{U}}L_{\mathcal{M}}\left\|\mu-\nu\right\|_{\mathrm{TV}}^{(\Gamma)} \\
&\le \frac{1}{\underline{\beta}\min_{\underline{\beta}\le\gamma\le\bar{\beta}}\mathbb{V}_\gamma\left[U(z,y,\mu)\right]}\left(\frac{\varepsilon\bar{U}}{\omega}+1\right)2\bar{\beta}e^{2\bar{\beta}\bar{U}}L_{\mathcal{M}}\left\|\mu-\nu\right\|_{\mathrm{TV}}^{(\Gamma)} \\
&\le \frac{\left(\frac{\varepsilon\bar{U}}{\omega}+1\right)2\bar{\beta}e^{\bar{\beta}2\bar{U}}L_{\mathcal{M}}}{\underline{\beta}\left(\min_{\underline{\beta}\le\gamma\le\bar{\beta}}\mathbb{V}_\gamma\left[U(z,y,\mu)\right]\wedge\min_{\underline{\beta}\le\gamma\le\bar{\beta}}\mathbb{V}_\gamma\left[U(z,y,\nu)\right]\right)}\left\|\mu-\nu\right\|_{\mathrm{TV}}^{(\Gamma)}.
\end{aligned}
\tag{48}
$$

We show that the variances appearing in the denominator (48) are strictly positive irrespective of $y$ and $\mu,\nu$, with which we can complete the proof.

Consider the family $\mathcal{Q}$ of Lipschitz continuous functions $f:D\to[0,+\infty)$ satisfying the following constraints:

**(a) Global bound**: Each $f \in Q$ has a global maximum of $b > 0$ and a global minimum of $0$.

**(b) Finiteness**: The set of points where $f$ achieves its global maximum, which is denoted by $M = \{x \in D | f(x) = b\}$, consists of a finite number of isolated points.

**(c) Local strict concavity:** $f$ is strictly concave in a neighborhood around each maximum $x' \in M$.

Fr is a fixed constant $a > 0$, and consider the following PDF: $p(x) = e^{af(x)} / \int_D e^{af(z)} \mathrm{d}z$, $f \in Q$. We prove the following **claim**: *the variance* $V_p[f]$ *of* $f$ *based on this PDF is strictly positive, i.e.,* $\inf_{f \in Q} V_p[f] > 0$.

This **claim** is proven as follows. Assume that $\inf_{f \in Q} V_p[f] = 0$ to lead to a contradiction. In this case, the random variable $f(X)$ (here, $X$ is valued in $D$) must be a constant, e.g., $\bar{f}$. If $V_p[f] \to 0$, then the push-forward measure defined by $p \circ f^{-1}$ must weakly converge to a Dirac delta $\delta_b$ concentrated at $b$. Because the PDF $p(x) \propto e^{af(x)}$ exponentially penalizes submaximum values of $f$, we need to have the convergence $f(X) \to b$ in the sense of probability; hence, we must have $\bar{f} = b$. This concentration toward $b$ requires the probability mass of the $\vartheta$-sublevel set, denoted by $S_\vartheta = \{x \in D | f(x) < b - \vartheta\}$, to vanish when $\vartheta \to 0$. However, because of the local strict concavity and finiteness of the maxima of $f$, the Lebesgue measure of $D - S_\vartheta$ shrinks to zero at a polynomial rate, $O(\vartheta^{n/2})$, with $n \in \mathbb{N}$ being the dimension of $D$, as $\vartheta \to 0$. Moreover, the contribution around the maxima vanishes: $\int_{D-S_\vartheta} e^{af(z)} \mathrm{d}z \to 0$ as $\vartheta \to 0$. Conversely, $S_\vartheta$ retains a strictly positive volume as $\vartheta \to 0$ because of the compactness of $D$, and its order is at $O(1)$. Consequently, the total integral over $S_\vartheta$ cannot approach 0. Therefore, $p$ is forced to maintain a nonzero spatial spread across regions where $f < b$, preventing $f(X)$ from collapsing into a constant, which contradicts $\inf_{f \in Q} V_p[f] = 0$. This implies that $\inf_{f \in Q} V_p[f] > 0$.

Now, **claim** stated above shows that there is a constant $I > 0$ independent of $y, \mu, \nu$ with

$$\min_{\underline{\beta} \le \gamma \le \bar{\beta}} \mathbb{V}_\gamma [U(z, y, \mu)] \wedge \min_{\underline{\beta} \le \gamma \le \bar{\beta}} \mathbb{V}_\gamma [U(z, y, \nu)] \ge I, \quad (49)$$

showing the Lipschitz continuity (23).

□

***Proof of Lemma 2***

The result follows due to the following estimate:

$$
\begin{aligned}
&\left|\mathbb{L}_{\beta(y,\mu)}(x,y,\mu)-\mathbb{L}_{\beta(y,\nu)}(x,y,\nu)\right| \\
&\le\left|\mathbb{L}_{\beta(y,\mu)}(x,y,\mu)-\mathbb{L}_{\beta(y,\nu)}(x,y,\mu)\right|+\left|\mathbb{L}_{\beta(y,\nu)}(x,y,\mu)-\mathbb{L}_{\beta(y,\nu)}(x,y,\nu)\right|
\end{aligned}, \tag{50}
$$

and we have

$$
\begin{aligned}
\left|\mathbb{L}_{\beta(y,\mu)}(x,y,\mu)-\mathbb{L}_{\beta(y,\nu)}(x,y,\mu)\right| &= \left|\frac{e^{\beta(y,\mu)U(x,y,\mu)}}{\int_D e^{\beta(y,\mu)U(m,y,\mu)}\mathrm{d}m}-\frac{e^{\beta(y,\nu)U(x,y,\mu)}}{\int_D e^{\beta(y,\nu)U(m,y,\mu)}\mathrm{d}m}\right| \\
&\le \left|\frac{e^{\beta(y,\mu)U(x,y,\mu)}}{\int_D e^{\beta(y,\mu)U(m,y,\mu)}\mathrm{d}m}-\frac{e^{\beta(y,\nu)U(x,y,\mu)}}{\int_D e^{\beta(y,\mu)U(m,y,\mu)}\mathrm{d}m}\right| \\
&\quad +\left|\frac{e^{\beta(y,\nu)U(x,y,\mu)}}{\int_D e^{\beta(y,\mu)U(m,y,\mu)}\mathrm{d}m}-\frac{e^{\beta(y,\nu)U(x,y,\mu)}}{\int_D e^{\beta(y,\nu)U(m,y,\mu)}\mathrm{d}m}\right| \\
&\le \left|e^{\beta(y,\mu)U(x,y,\mu)}-e^{\beta(y,\nu)U(x,y,\mu)}\right| \\
&\quad +e^{\beta(y,\nu)U(x,y,\mu)}\int_D\left|e^{\beta(y,\mu)U(m,y,\mu)}-e^{\beta(y,\nu)U(m,y,\mu)}\right|\mathrm{d}m \\
&\le \bar{U}e^{\bar{\beta}\bar{U}}\left(1+e^{\bar{\beta}\bar{U}}\right)\|\mu-\nu\|_{\mathrm{TV}}^{(\Gamma)} \\
&:= C_{0,1}\|\mu-\nu\|_{\mathrm{TV}}^{(\Gamma)}
\end{aligned} \tag{51}
$$

and similarly, there is a constant $C_{0,2}$ depending on $\bar{\beta},\bar{U},L_{\mathcal{M}}$ such that

$$
\left|\mathbb{L}_{\beta(y,\nu)}(x,y,\mu)-\mathbb{L}_{\beta(y,\nu)}(x,y,\nu)\right|\le C_{0,2}\|\mu-\nu\|_{\mathrm{TV}}^{(\Gamma)}. \tag{52}
$$

We complete the proof by setting $C_0=C_{0,1}+C_{0,2}$.

□

***Proof of Proposition 1***

The proof is an adaptation of Proof of Theorem 3.4 in Lahkar and Riedel [3] to our case. First, we consider the auxiliary, truncated logit dynamic: For any Borel measurable set $A$ in $D$:

$$
\frac{\mathrm{d}\mu_t(A,y)}{\mathrm{d}t}=\left(2-\|\mu_t\|_{\mathrm{TV}}^{(\Gamma)}\right)_+\left(\int_A\mathbb{L}_{\beta(y,\mu_t)}(x,y,\mu)\mathrm{d}z-\mu_t(A,y)\right),\quad t>0,\quad y\in\Gamma \tag{53}
$$

with the initial condition being unchanged. This equation has an artificial factor $\left(2-\|\mu_t\|_{\mathrm{TV}}^{(\Gamma)}\right)_+$ that equals 1 if the solution belongs to $\mathcal{P}^{(\Gamma)}$. With the help of the Lipschitz continuity result in **Lemma 2**, the right-hand side of (53) is Lipschitz continuous with respect to $\mu_t\in\mathcal{M}^{(\Gamma)}$ because the artificial factor vanishes when $\mu_t\in\mathcal{M}_2^{(\Gamma)}$. Applying the generalized Picard–Lindelöf theorem (Lahkar et al. [15]; Yoshioka [16]; see also p.79 in Zeidler [17]) to (53) shows that this equation admits a unique solution $(\mu_t)_{t\ge 0}$ such that $\mu_t\in\mathcal{M}^{(\Gamma)}$ for all $t\ge 0$. The fact that $\|\mu_t\|_{\mathrm{TV}}^{(\Gamma)}=1$ follows from integrating equation (53) in time, and then we obtain that the logit dynamic (19) admits a unique solution $\mu$ such that $\mu_t\in\mathcal{M}^{(\Gamma)}$ and

$\|\mu_t\|_{\mathrm{TV}}^{(\Gamma)} = 1$ for all $t \geq 0$, and then we can formally drop the artificial factor $\left(2 - \|\mu_t\|_{\mathrm{TV}}^{(\Gamma)}\right)_+$. Finally, the nonnegativity of this solution follows from the following variation of the constant formula whose right-hand side is nonnegative:

$$\mu_t(A,y) = \mu_0(A,y)e^{-t} + \int_0^t e^{-(t-s)} \int_A \mathbb{L}_{\beta(y,\mu_s)}(z,y,\mu)\mathrm{d}s\mathrm{d}z\,,\quad t > 0\,,\quad y \in \Gamma \tag{54}$$

for any Borel measurable set $A$ in $D$. Finally, the local Lipschitz continuity of the right-hand side of the logit dynamic (18) for any $\mu_t \in \mathcal{P}^{(\Gamma)}$ and $t \geq 0$ shows that the solution found in this proof is the desired unique solution to (18).

□

***Proof of Lemma 3***

The proof is essentially the same as those for **Lemmas 1-2**; here, we only need to formally replace the total variation norm with the Wasserstein distance.

□

***Proof of Proposition 2***

Assume that there are two strong solutions $\{X_t^y\}_{t\in[0,T],\, y\in\Gamma}$ and $\{Z_t^y\}_{t\in[0,T],\, y\in\Gamma}$ to (20) with the common initial condition $\{X_0^y\}_{y\in\Gamma}$. We temporally fix $y \in \Gamma$ and $t \in (0,T)$. In the proof, we omit the arguments of $\mathbb{L}_\beta$ and $\mathcal{L}$ when there is no confusion. Note that the boundedness of the solutions, i.e., $X_t^y, Z_t^y \in D$ is clear from the SDE form (20).

By assumption, we have

$$X_t^y = X_0^y + \int_0^t \int_D \int_0^\infty \left(z - X_s^y\right) \mathbb{I}_{\{w \leq \mathbb{L}_{\beta(y,\mathcal{L}_{s-}(X))}\}} N\left(\mathrm{d}s, \mathrm{d}z, \mathrm{d}w, \{y\}\right) \tag{55}$$

and

$$Z_t^y = X_0^y + \int_0^t \int_D \int_0^\infty \left(z - Z_s^y\right) \mathbb{I}_{\{w \leq \mathbb{L}_{\beta(y,\mathcal{L}_{s-}(Z))}\}} N\left(\mathrm{d}s, \mathrm{d}z, \mathrm{d}w, \{y\}\right). \tag{56}$$

Hence, obtain

$$\begin{aligned}
\left|X_t^y - Z_t^y\right| &\leq \int_0^t \int_D \int_0^\infty \left| \left(z - X_s^y\right) \mathbb{I}_{\{w \leq \mathbb{L}_{\beta(y,\mathcal{L}_{s-}(X))}\}} - \left(z - Z_s^y\right) \mathbb{I}_{\{w \leq \mathbb{L}_{\beta(y,\mathcal{L}_{s-}(Z))}\}} \right| N\left(\mathrm{d}s, \mathrm{d}z, \mathrm{d}w, \{y\}\right) \\
&\leq \int_0^t \int_D \int_0^\infty \left| \left(z - X_s^y\right) \mathbb{I}_{\{w \leq \mathbb{L}_{\beta(y,\mathcal{L}_{s-}(X))}\}} - \left(z - Z_s^y\right) \mathbb{I}_{\{w \leq \mathbb{L}_{\beta(y,\mathcal{L}_{s-}(X))}\}} \right| N\left(\mathrm{d}s, \mathrm{d}z, \mathrm{d}w, \{y\}\right) \\
&\quad + \int_0^t \int_D \int_0^\infty \left| \left(z - Z_s^y\right) \mathbb{I}_{\{w \leq \mathbb{L}_{\beta(y,\mathcal{L}_{s-}(X))}\}} - \left(z - Z_s^y\right) \mathbb{I}_{\{w \leq \mathbb{L}_{\beta(y,\mathcal{L}_{s-}(Z))}\}} \right| N\left(\mathrm{d}s, \mathrm{d}z, \mathrm{d}w, \{y\}\right) \\
&\leq \int_0^t \int_D \int_0^\infty \left|Z_s^y - X_s^y\right| \mathbb{I}_{\{w \leq \mathbb{L}_{\beta(y,\mathcal{L}_{s-}(X))}\}} N\left(\mathrm{d}s, \mathrm{d}z, \mathrm{d}w, \{y\}\right) \\
&\quad + \int_0^t \int_D \int_0^\infty \left|z - Z_s^y\right| \left| \mathbb{I}_{\{w \leq \mathbb{L}_{\beta(y,\mathcal{L}_{s-}(X))}\}} - \mathbb{I}_{\{w \leq \mathbb{L}_{\beta(y,\mathcal{L}_{s-}(Z))}\}} \right| N\left(\mathrm{d}s, \mathrm{d}z, \mathrm{d}w, \{y\}\right) \\
&:= I_1 + I_2.
\end{aligned} \tag{57}$$

We evaluate the last two terms in (57). The first term $I_1$ is evaluated as follows (exchange between integration and expectation is allowed because of the strict boundedness of the integrant):

$$\begin{aligned}\mathbb{E}[I_1] &= \int_0^t\int_D \mathbb{E}\left[\left|Z_s^y - X_s^y\right|\mathbb{L}_{\beta(y,\mathcal{L}_{s-}(X))}\left(z,y,\mathcal{L}_s(X)\right)\right]\mathrm{d}s\mathrm{d}z \\ &\leq \exp\left(\bar{\beta}\bar{U}\right)\int_0^t\int_D \mathbb{E}\left[\left|Z_s^y - X_s^y\right|\right]\mathrm{d}s\mathrm{d}z \\ &= \exp\left(\bar{\beta}\bar{U}\right)\int_0^t \mathbb{E}\left[\left|Z_s^y - X_s^y\right|\right]\mathrm{d}s.\end{aligned} \tag{58}$$

For the second term $I_2$, we have

$$\begin{aligned}I_2 &= \int_0^t\int_D\int_0^\infty \left|z - Z_s^y\right|\left|\mathbb{I}_{\{w\leq\mathbb{L}_{\beta(y,\mathcal{L}_{s-}(X))}\}} - \mathbb{I}_{\{w\leq\mathbb{L}_{\beta(y,\mathcal{L}_{s-}(Z))}\}}\right| N\left(\mathrm{d}s,\mathrm{d}z,\mathrm{d}w,\{y\}\right) \\ &\leq \int_0^t\int_D\int_0^\infty \left|\mathbb{I}_{\{w\leq\mathbb{L}_{\beta(y,\mathcal{L}_{s-}(X))}\}} - \mathbb{I}_{\{w\leq\mathbb{L}_{\beta(y,\mathcal{L}_{s-}(Z))}\}}\right| N\left(\mathrm{d}s,\mathrm{d}z,\mathrm{d}w,\{y\}\right)\end{aligned} \tag{59}$$

because of the assumption about $D$. Afterward, we proceed as follows:

$$\begin{aligned}\mathbb{E}[I_2] &\leq \int_0^t\int_D \mathbb{E}\left[\int_0^{+\infty}\left|\mathbb{I}_{\{w\leq\mathbb{L}_{\beta(y,\mathcal{L}_s(X))}\}} - \mathbb{I}_{\{w\leq\mathbb{L}_{\beta(y,\mathcal{L}_s(Z))}\}}\right|\mathrm{d}w\right]\mathrm{d}s\mathrm{d}z \\ &= \int_0^t\int_D \mathbb{E}\left[\mathbb{L}_{\beta(y,\mathcal{L}_s(X))}\left(z,y,\mathcal{L}_s(X)\right) - \mathbb{L}_{\beta(y,\mathcal{L}_s(Z))}\left(z,y,\mathcal{L}_s(Z)\right)\right]\mathrm{d}s\mathrm{d}z\end{aligned}. \tag{60}$$

Here, we used the following elementary identity: for any $a,b>0$,

$$\int_0^{+\infty}\left|\mathbb{I}_{\{w\leq a\}} - \mathbb{I}_{\{w\leq b\}}\right|\mathrm{d}w = |a-b|. \tag{61}$$

Then, the Lipschitz estimate of the logit function $\mathbb{L}_\beta$ (**Lemma 3**) shows that there exists a global constant $C_1>0$ independent of the solutions, $z\in D$, and $T$, such that

$$\left|\mathbb{L}_{\beta(y,\mathcal{L}_s(X))}\left(z,y,\mathcal{L}_s(X)\right) - \mathbb{L}_{\beta(y,\mathcal{L}_s(Z))}\left(z,y,\mathcal{L}_s(Z)\right)\right| \leq C_1 W_1\left(\mathcal{L}\left(X_s^y\right),\mathcal{L}\left(Z_s^y\right)\right). \tag{62}$$

With the help of the inequality (9), we obtain

$$\begin{aligned}\mathbb{E}[I_2] &\leq C_1\int_0^t\int_D \mathbb{E}\left[W_1\left(\mathcal{L}\left(X_s^y\right),\mathcal{L}\left(Z_s^y\right)\right)\right]\mathrm{d}s\mathrm{d}z \\ &\leq C_1\int_0^t\int_D \mathbb{E}\left[\left|X_s^y - Z_s^y\right|\right]\mathrm{d}s\mathrm{d}z \\ &= C_1\int_0^t \mathbb{E}\left[\left|X_s^y - Z_s^y\right|\right]\mathrm{d}s.\end{aligned} \tag{63}$$

In summary, we have

$$\mathbb{E}\left[\left|X_t^y - Z_t^y\right|\right] \leq C_1\exp\left(\bar{\beta}\bar{U}\right)\int_0^t \mathbb{E}\left[\left|X_s^y - Z_s^y\right|\right]\mathrm{d}s \tag{64}$$

and hence

$$\begin{aligned}\sup_{0\leq r\leq t}\mathbb{E}\left[\left|X_r^y - Z_r^y\right|\right] &\leq C_1\exp\left(\bar{\beta}\bar{U}\right)\sup_{0\leq r\leq t}\int_0^r \mathbb{E}\left[\left|X_s^y - Z_s^y\right|\right]\mathrm{d}s \\ &\leq C_1\exp\left(\bar{\beta}\bar{U}\right)\int_0^t \sup_{0\leq r\leq s}\mathbb{E}\left[\left|X_r^y - Z_r^y\right|\right]\mathrm{d}s\end{aligned} \tag{65}$$

Integrating both sides of (64) for $y\in\Gamma$ and applying the classical Gronwall inequality yields

$$\int_\Gamma \sup_{0\leq r\leq t}\mathbb{E}\left[\left|X_r^y - Z_r^y\right|\right]\mathrm{d}y = 0,\quad 0\leq t\leq T, \tag{66}$$

which proves the proposition because $T>0$ is arbitrary.

$\square$

***Proof of Proposition 3***

Fix $T>0$. The proof is based on a Picard iteration method. Set some $X^{(0)}\in\mathcal{H}_T$ such that $X_0^{(0),y}=X_0^y\in D$ ($y\in\Gamma$). Afterward, we recursively set the SDEs: for $0<t<T$,

$$X_t^{(n+1),y}=X_0^y+\int_0^t\int_D\int_0^\infty\left(z-X_s^{(n+1),y}\right)\mathbb{I}_{\left\{w\le\mathbb{L}_{\beta\left(y,\mathcal{L}_{s-}\left(X^{(n)}\right)\right)}\right\}}N\left(\mathrm{d}s,\mathrm{d}z,\mathrm{d}w,\{y\}\right),\quad n=0,1,2,\ldots. \tag{67}$$

For each $n=0,1,2,...$, this is a jump-driven SDE driven without any mean-field effects; hence, each $X^{(n+1)}\in\mathcal{H}_T$ exists uniquely. Note that the boundedness of the solutions, i.e., $X_t^{(n+1),y}\in D$, is clear from (67) because $X_t^{(n+1),y}$ is generated by $\mathbb{L}_{\beta\left(y,\mathcal{L}_{t-}\left(X^{(n)}\right)\right)}$.

We consider the SDE (67) for some $n$ and $n-1$ ($n\ge1$) to obtain

$$\begin{aligned}&\mathbb{E}\left[\left|X_t^{(n+1),y}-X_t^{(n),y}\right|\right]\\&\le\mathbb{E}\left[\int_0^t\int_D\int_0^\infty\left|\left(z-X_s^{(n+1),y}\right)\mathbb{I}_{\left\{w\le\mathbb{L}_{\beta\left(y,\mathcal{L}_{s-}\left(X^{(n)}\right)\right)}\right\}}-\left(z-X_s^{(n),y}\right)\mathbb{I}_{\left\{w\le\mathbb{L}_{\beta\left(y,\mathcal{L}_{s-}\left(X^{(n-1)}\right)\right)}\right\}}\right|N\left(\mathrm{d}s,\mathrm{d}z,\mathrm{d}w,\{y\}\right)\right].\end{aligned} \tag{68}$$

By an argument analogous to that used in **Proof of Proposition 2**, we arrive at the following inequality:

$$\mathbb{E}\left[\left|X_t^{(n+1),y}-X_t^{(n),y}\right|\right]\le\exp\left(\bar{\beta}\bar{U}\right)\int_0^t\mathbb{E}\left[\left|X_s^{(n+1),y}-X_s^{(n),y}\right|\right]\mathrm{d}s+C_1\int_0^t\mathbb{E}\left[\left|X_s^{(n),y}-X_s^{(n-1),y}\right|\right]\mathrm{d}s. \tag{69}$$

Now, we set

$$\Lambda^{(n)}(t)=\int_\Gamma\sup_{0\le r\le t}\mathbb{E}\left[\left|X_r^{(n),y}-X_r^{(n-1),y}\right|\right]\mathrm{d}y,\quad n=1,2,3,\ldots. \tag{70}$$

Then, by (69), we have

$$\begin{aligned}&\sup_{0\le r\le t}\mathbb{E}\left[\left|X_r^{(n+1),y}-X_r^{(n),y}\right|\right]\\&\le\exp\left(\bar{\beta}\bar{U}\right)\sup_{0\le r\le t}\int_0^r\mathbb{E}\left[\left|X_r^{(n+1),y}-X_r^{(n),y}\right|\right]\mathrm{d}s+C_1\sup_{0\le r\le t}\int_0^r\mathbb{E}\left[\left|X_r^{(n),y}-X_r^{(n-1),y}\right|\right]\mathrm{d}s\\&\le\exp\left(\bar{\beta}\bar{U}\right)\int_0^t\sup_{0\le r\le s}\mathbb{E}\left[\left|X_r^{(n+1),y}-X_r^{(n),y}\right|\right]\mathrm{d}s+C_1\int_0^t\sup_{0\le r\le s}\mathbb{E}\left[\left|X_r^{(n),y}-X_r^{(n-1),y}\right|\right]\mathrm{d}s,\end{aligned} \tag{71}$$

and hence

$$\Lambda^{(n+1)}(t)\le\exp\left(\bar{\beta}\bar{U}\right)\int_0^t\Lambda^{(n+1)}(s)\mathrm{d}s+C_1\int_0^t\Lambda^{(n)}(s)\mathrm{d}s. \tag{72}$$

Because each $\Lambda^{(n)}(t)$ is increasing for $t\in[0,T]$, by the classical Gronwall inequality, we obtain

$$\begin{aligned}\Lambda^{(n+1)}(t)&\le\exp\left(\int_0^t\exp\left(\bar{\beta}\bar{U}\right)\mathrm{d}s\right)C_1\int_0^t\Lambda^{(n)}(s)\mathrm{d}s\\&\le C_1\exp\left(T\exp\left(\bar{\beta}\bar{U}\right)\right)\int_0^t\Lambda^{(n)}(s)\mathrm{d}s\\&\left(:=C_2\int_0^t\Lambda^{(n)}(s)\mathrm{d}s\right).\end{aligned} \tag{73}$$

Then, owing to the standard iterative argument (e.g., Proof of Theorem 4.6.7 in Méléard [65]), with $C_3 = \int_\Gamma \sup_{0 \le r \le T} \mathbb{E}\left[\left|X_r^{(1),y} - X_r^{(0),y}\right|\right] \mathrm{d}y < +\infty$ , we have

$$\Lambda^{(n)}(t) \le C_3 \frac{(C_2 T)^{n-1}}{(n-1)!}, \quad n = 1,2,3,... \tag{74}$$

and hence we obtain

$$\sum_{n=0}^{+\infty} \Lambda^{(n)}(t) \le \sum_{n=0}^{+\infty} C_3 \frac{(C_2 T)^n}{(n)!} \le C_3 \exp(C_2 T) < +\infty . \tag{75}$$

This inequality shows that $X^{(n)}$ ($n = 0,1,2,...$) is a Cauchy sequence in $\mathcal{H}_T$; hence, there exists a limit $X^{(\infty)} \in \mathcal{H}_T$ such that

$$\lim_{n \to +\infty} \left\| X^{(n)} - X^{(\infty)} \right\|_{\mathcal{H}_T} = 0 . \tag{76}$$

Therefore, we also have the following weak convergence:

$$\lim_{n \to +\infty} \sup_{0 \le t \le T} W_1\left(\mathcal{L}\left(X_t^{(n),y}\right), \mathcal{L}\left(X_t^{(\infty),y}\right)\right) = 0 \quad \text{at almost all} \quad y \in \Gamma , \tag{77}$$

and hence obtain the convergence of the indicator functions:

$$\begin{aligned}
&\int_\Gamma \sup_{0 \le t \le T} \mathbb{E}\left[\int_0^\infty \left| \mathbb{I}_{\{w \le \mathbb{L}_{\beta(y,\mathcal{L}_t(X^{(n)}))}\}} - \mathbb{I}_{\{w \le \mathbb{L}_{\beta(y,\mathcal{L}_t(X^{(\infty)}))}\}} \right|\right] \mathrm{d}y \\
&= \int_\Gamma \sup_{0 \le t \le T} \mathbb{E}\left[\left| \mathbb{L}_{\beta(y,\mathcal{L}_t(X^{(n)}))}\left(X_t^{(n),y}\right) - \mathbb{L}_{\beta(y,\mathcal{L}_t(X^{(\infty)}))}\left(X_t^{(\infty),y}\right) \right|\right] \mathrm{d}y \\
&\le C_1 \int_\Gamma \sup_{0 \le t \le T} \mathbb{E}\left[W_1\left(\mathcal{L}\left(X_t^{(n),y}\right), \mathcal{L}\left(X_t^{(\infty),y}\right)\right)\right] \mathrm{d}y \\
&\le C_1 \int_\Gamma \mathbb{E}\left[\sup_{0 \le t \le T} W_1\left(\mathcal{L}\left(X_t^{(n),y}\right), \mathcal{L}\left(X_t^{(\infty),y}\right)\right)\right] \mathrm{d}y \\
&\to 0.
\end{aligned} \tag{78}$$

These observations imply that the limit $X^{(\infty)} \in \mathcal{H}_T$ solves the SDE (20) and hence completes the proof.

□

***Proof of Proposition 4.***

The proof is the direct application of the classical Itô's formula (Lemma 5 in Kuehn and Pulido [37]), showing that $\left\{\mathcal{L}_t\left(X^\cdot\right)\right\}_{t \ge 0}$ satisfies the logit dynamic (18); in our case, applying Itô's formula for jump processes suffices (e.g., Theorem 1.16 in Øksendal and Sulem [34]).

□

## B. Numerical method

### B.1 Discretization scheme

A naïve finite difference scheme to discretize our logit dynamic is presented, where we assume $D=[0,1]$, but its multidimensional and graph extensions are straightforward because the dynamic does not involve spatial derivatives. The scheme is that in Appendix of Yoshioka [26] with an adaptation to our setting.

The time increment is set as $\Delta t>0$, and the spatial increment as $\Delta x=\Delta y=1/N$ with some $N\in\mathbb{N}$. We prepare one-dimensional cells $d_i=[(i-1)1/N,i/N)$ for $i=1,2,...,N-1$ and $d_N=[1-1/N,1]$. We also set two-dimensional cells $D_{i,j}=d_i\times d_j$ ($i,j=1,2,...,N$). The quantity discretized at time $k\Delta t$ ($k=0,1,2,...$) and in $D_{i,j}$ ($i,j=1,2,...,N$) is denoted by the subscript $i,j,k$. The initial condition $\mu_0$ in $D_{i,j}$ is $\mu_{0,i,j}=\int_{D_{i,j}}\mu_0(\mathrm{d}x)$. Utilities $U$ are approximated by using a common midpoint rule evaluated at the centers of $D_{i,j}$.

Each $\mu_{i,j,k}$ ($i,j=1,2,...,N$) is computed as follows:

$$\mu_{i,j,k}=(1-\Delta t)\mu_{i,j,k-1}+\Delta t\frac{\exp\left(\beta_{j,k-1}U_{i,j,k-1}\right)}{\sum_{l=1}^{N}\exp\left(\beta_{j,k-1}U_{l,j,k-1}\right)\Delta x}, \tag{79}$$

where $\beta_{j,k-1}$ is uniquely found from the following equation:

$$0=\sum_{i=1}^{N}\left\{\frac{\exp\left(\beta_{j,k-1}U_{i,j,k-1}\right)}{\sum_{l=1}^{N}\exp\left(\beta_{j,k-1}U_{l,j,k-1}\right)\Delta x}\ln\left(\frac{\exp\left(\beta_{j,k-1}U_{i,j,k-1}\right)}{\sum_{l=1}^{N}\exp\left(\beta_{j,k-1}U_{l,j,k-1}\right)\Delta x}\right)\right\}\Delta x+\omega\beta_{j,k-1}-\varepsilon\ \left(:=g\left(\beta_{j,k-1}\right)\right). \tag{80}$$

After an elementary calculation, we have

$$\frac{\mathrm{d}g(\beta)}{\mathrm{d}\beta}=\sum_{i=1}^{N}\left\{\frac{\left(U_{i,j,k-1}\right)^2\exp\left(\beta U_{i,j,k-1}\right)}{\sum_{l=1}^{N}\exp\left(\beta U_{l,j,k-1}\right)\Delta x}-\left(\frac{U_{i,j,k-1}\exp\left(\beta U_{i,j,k-1}\right)}{\sum_{l=1}^{N}\exp\left(\beta U_{l,j,k-1}\right)\Delta x}\right)^2\right\}\Delta x+\omega\geq\omega>0,\quad \beta>0, \tag{81}$$

where the summation term in (80) represents a discretization of the variance of a utility against a logit function (i.e., the discrete version of $\mathbb{V}_\beta\left[U(z,y,\mu)\right]$ in **Proof of Lemma 1**). The strict positivity result (81) implies that the following Newton iteration never breaks down by the sign change of the denominator:

$$\beta_{j,k-1}\to\beta_{j,k-1}-\left(\frac{\mathrm{d}g\left(\beta_{j,k-1}\right)}{\mathrm{d}\beta}\right)^{-1}g\left(\beta_{j,k-1}\right)\ \left(:=J\left(\beta_{j,k-1}\right)\right). \tag{82}$$

This is another theoretical advantage of incorporating the regularization term into the cost constraint. The system (80) is iteratively solved for $\eta_{k-1}$ by using **Algorithm 1**. For all $j$, the initial guess $\beta^{(0)}_{j,k-1}$ is 10 and is updated recursively for $k>0$ by using the converged value at the previous time step.

***Algorithm 1***

***Step 0.*** *Set the error threshold* $\mathrm{Er}\left(=10^{-12}\right)$/

***Step 1.*** *Set* $\beta_{j,k-1}^{(0)} > 0$. *Set* $n = 0$.

***Step 2.*** *Use the Newton iteration* $\beta_{j,k-1}^{(n+1)} = J\left(\beta_{j,k-1}^{(n)}\right)$.

***Step 3.*** *Compute the error* $\mathrm{Er}^{(n)} = \max_{j}\left|g\left(\beta_{j,k-1}\right)\right|$.

***Step 4.*** *If* $\mathrm{Er}^{(n)} \le \mathrm{Er}$, *then output* $\beta_{j,k-1}^{(n+1)}$ *as the numerical solution to* (80) *and terminate the algorithm. If* $\mathrm{Er}^{(n)} > \mathrm{Er}$, *then set* $n \to n+1$ *and go to* ***Step 2***.

**B.2 Convergence study**

We conduct a convergence study of the finite difference method to investigate its convergence speed and to show that the computational resolutions $\Delta x = \Delta y = 1/400$ and $\Delta t = 1/200$ are sufficiently fine. Here, we use the nominal parameter values specified in **Section 4.1**. Because the equilibrium $\hat{x}_y$ of (26) is available explicitly and is closely related to the action average $\bar{x}_y$ defined in (29), we compare the action averages for different computational resolutions and compare them at a stationary state. Typically, fewer than ten Newton iterations suffice at each time step and each $j$ for the cases studied in this paper.

Here, the "reference" implies the action average computed using the resolution $\Delta x = \Delta y = 1/1600$ and $\Delta t = 1/200$. **Table 1** compares the maximum absolute difference between the reference and numerical $\bar{x}_y$ at a stationary state for different computational resolutions, demonstrating that the convergence speed is second order in $N$ (the absolute difference is almost halved by doubling $N$). Moreover, the difference is considered sufficiently small at $N = 400$ used in the main text. We also examine the difference in time by comparing the maximum and minimum of the computed PDFs at each time step. We examine different values of $\Delta t$ with $N$ fixed. Here, the "reference" implies the action average computed using the resolution $\Delta x = \Delta y = 1/400$ and $\Delta t = 1/800$. **Table 2** compares the maximum absolute difference between the reference and numerical results during the time interval $\left(0,15\right]$, after which the numerical solutions are almost stationary. The computational results suggest that the convergence speed is approximately first order in $1/\Delta t$, and the use of $\Delta t = 1/200$ suffices for the analysis in the main text.

Finally, for the second application in the main text, **Figure B1** compares the computed PDFs of the arrival intensity at the beginning on July 1 with the resolutions $\Delta x = \Delta y = 1/400$ with $\Delta t = 1/24$ and $\Delta x = \Delta y = 1/800$ with $\Delta t = 1/48$. The computational results are not distinguishable from each other, and their difference is 0.13 in the maximum norm, which amounts to a couple of percentages in relative error. This error does not critically affect the results obtained in the second application.

**Table 1.** Absolute difference between the reference and other numerical $\bar{x}_y$ at a stationary state for different $N$.

| $N$ | Absolute difference |
|---|---|
| 100 | 9.231E-07 |
| 200 | 2.292E-07 |
| 400 | 5.471E-08 |
| 800 | 1.096E-08 |

**Table 2.** Absolute difference between the reference and numerical $\mathrm{AC}_k$ for different $\Delta t$.

| | Absolute difference | |
|---|---|---|
| $\Delta t$ | Maximum PDF | Minimum PDF |
| 1/50 | 1.192.E-02 | 3.480.E-03 |
| 1/100 | 5.541.E-03 | 1.617.E-03 |
| 1/200 | 2.370.E-03 | 6.916.E-04 |
| 1/400 | 7.891.E-04 | 2.303.E-04 |

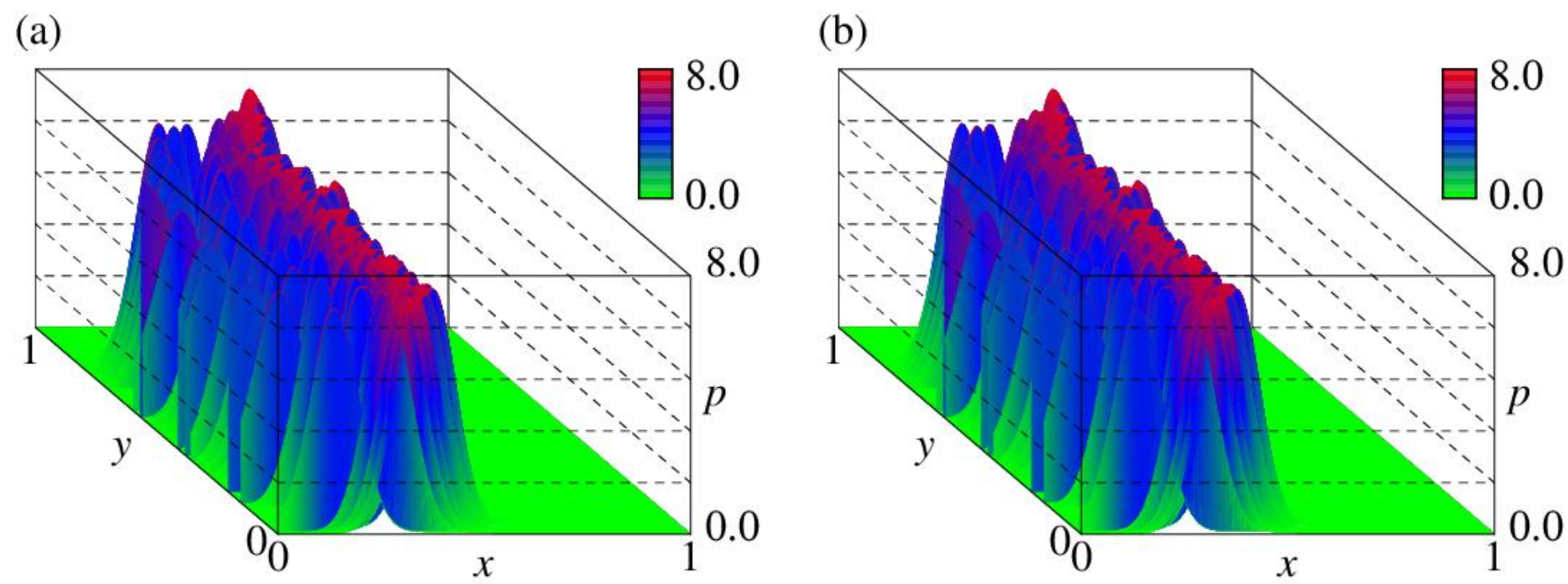

**Figure B1.** Computed PDFs of the arrival intensity at the beginnings on July 1 with (a) $\Delta x = \Delta y = 1/400$ and $\Delta t = 1/24$ (b) $\Delta x = \Delta y = 1/800$ and $\Delta t = 1/48$.